\documentclass[11pt]{article}
\usepackage{amssymb}
\usepackage{amsmath}
\makeatletter
\@addtoreset{equation}{section}
\makeatother

\title{Quenched Large Deviations for one dimensional
Nonlinear Filtering}

\author{\'Etienne Pardoux
\thanks{Member of the Institut Universitaire de France}\\
LATP, Univ. de Provence and CNRS\\
CMI, 39 rue Joliot Curie\\
13453 Marseille Cedex 13, France\\
email: pardoux@cmi.univ-mrs.fr
\and
 Ofer Zeitouni\thanks{Part of this work was done while visiting the 
LATP, University of Provence and CNRS, Marseille. It was also
supported by  
NSF grant number
DMS-0302230.}\\
Dept.
 of Electrical Engineering,
Technion\\
Haifa 32000, Israel\\ and\\
Dept. of Mathematics, 
Vincent Hall\\
University of Minnesota\\
Minneapolis 55455, USA\\
email: zeitouni@math.umn.edu}
\date{May 21, 2003. Revised November 21, 2003 and May 30, 2004.}
\begin{document}

% french
\newtheorem{theorem*}     {theorem}
\newtheorem{proposition}  [theorem]{Proposition}
\newtheorem{definition}   [theorem]{D\'efinition}
\newtheorem{lemme}        [theorem]{Lemme}
\newtheorem{corollaire}   [theorem]{Corollaire}
\newtheorem{resultat}     [theorem]{R\'esultat}
\newtheorem{eexercice}    [theorem]{Exercice}
\newtheorem{rremarque}    [theorem]{Remarque}
\newtheorem{pprobleme}    [theorem]{Probl\`eme}
\newtheorem{eexemple}     [theorem]{Exemple}
\newcommand{\preuve}      {\paragraph{Preuve}}
\newenvironment{probleme} {\begin{pprobleme}\rm}{\end{pprobleme}}
\newenvironment{remarque} {\begin{rremarque}\rm}{\end{rremarque}}
\newenvironment{exercice} {\begin{eexercice}\rm}{\end{eexercice}}
\newenvironment{exemple}  {\begin{eexemple}\rm}{\end{eexemple}}
%
% english
%
\newtheorem{e-theorem}      {Theorem}[section]
\newtheorem{e-theorem*}     {Theorem}
\newtheorem{e-resultat*}    [e-theorem*]{Result}
\newtheorem{e-proposition}  [theorem]{Proposition}
\newtheorem{e-definition}   [theorem]{Definition}
\newtheorem{e-lemme}        [theorem]{Lemma}
\newtheorem{e-corollaire}   [theorem]{Corollary}
\newtheorem{e-resultat}     [theorem]{Result}
\newtheorem{e-eexercice}    [theorem]{Exercise}
\newtheorem{e-rremarque}    [theorem]{Remark}
\newtheorem{e-pprobleme}    [theorem]{Problem}
\newtheorem{e-eexemple}     [theorem]{Example}
\newcommand{\proof}         {\paragraph{Proof~: }}
\newcommand{\hint}          {\paragraph{Hint}}
\newcommand{\heuristicproof}{\paragraph{heuristic proof}}
\newenvironment{e-probleme} {\begin{e-pprobleme}\rm}{\end{e-pprobleme}}
\newenvironment{e-remarque} {\begin{e-rremarque}\rm}{\end{e-rremarque}}
\newenvironment{e-exercice} {\begin{e-eexercice}\rm}{\end{e-eexercice}}
\newenvironment{e-exemple}  {\begin{e-eexemple}\rm}{\end{e-eexemple}}
\newcommand{\reell}    {{{\rm I\! R}^l}}
\newcommand{\reels}    {{{\rm I\! R}}}
\newcommand{\rationels}    {{{\rm I\!\!\! Q}}}
\newcommand{\reeln}    {{{\rm I\! R}^n}}
\newcommand{\reelk}    {{{\rm I\! R}^k}}
\newcommand{\reelm}    {{{\rm I\! R}^m}}
\newcommand{\reelp}    {{{\rm I\! R}^p}}
\newcommand{\reeld}    {{{\rm I\! R}^d}}
\newcommand{\reeldd}   {{{\rm I\! R}^{d\times d}}}
\newcommand{\reelnn}   {{{\rm I\! R}^{n\times n}}}
\newcommand{\reelnd}   {{{\rm I\! R}^{n\times d}}}
\newcommand{\reeldn}   {{{\rm I\! R}^{d\times n}}}
\newcommand{\reelkd}   {{{\rm I\! R}^{k\times d}}}
\newcommand{\reelkl}   {{{\rm I\! R}^{k\times l}}}
\newcommand{\reelN}    {{{\rm I\! R}^N}}
\newcommand{\reelM}    {{{\rm I\! R}^M}}
\newcommand{\reelplus} {{{\rm I\! R}^+}}
\newcommand{\reelo}    {{{\rm I\! R}\setminus\{0\}}}
\newcommand{\reld}    {{{\rm I\! R}_d}}
\newcommand{\relplus} {{{\rm I\! R}_+}}
\newcommand{\1}        {{\bf 1}}

\newcommand{\cov}      {{\hbox{cov}}}
\newcommand{\sss}      {{\cal S}}
\newcommand{\indic}    {{{\rm I\!\! I}}}
\newcommand{\pp}       {{{\rm I\! P}}}
\newcommand{\qq}       {{{\rm I\!\!\! Q}}}
\newcommand{\ee}       {{{\rm I\! E}}}

\newcommand{\B}        {{\hbox{\rm I\hskip -0.2cm B}}}
\newcommand{\cc}       {{\hbox{\rm I$\!\!\!$ C}}}
\newcommand{\HHH}       {{\hbox{\rm I\hskip -0.2cm H}}}
\newcommand{\N}        {{  {\rm I\hskip -2pt N}}}
\newcommand{\R}        {{ {\rm I \hskip -2pt R}}}
\newcommand{\D}        {{ {\rm I \hskip -2pt D}}}
\newcommand{\Q}       {{{\rm I\!\!\! Q}}}
\newcommand{\C}        {{\bf C}}	% ensemble des nombres complexes
\newcommand{\T}        {{\bf T}}	% espace de la variable temporelle
\newcommand{\E}        {{ {\rm I \hskip -2pt E}}}	% esp math
\renewcommand{\P}        {{ {\rm I \hskip -2pt P}}}
\newcommand{\h}   {\overline{h}}
\newcommand{\rfr}[1]    {\stackrel{\circ}{#1}}
\newcommand{\equiva}    {\displaystyle\mathop{\simeq}}
\newcommand{\eqdef}     {\stackrel{\triangle}{=}}
\newcommand{\limps}     {\mathop{\hbox{\rm lim--p.s.}}}
\newcommand{\Limsup}    {\mathop{\overline{\rm lim}}}
\newcommand{\Liminf}    {\mathop{\underline{\rm lim}}}
\newcommand{\Inf}       {\mathop{\rm Inf}}
%                       \joinrel\rightarrow}\;}}
\newcommand{\vers}      {\mathop{\;{\rightarrow}\;}}
\newcommand{\versup}    {\mathop{\;{\nearrow}\;}}
\newcommand{\versdown}  {\mathop{\;{\searrow}\;}}
\newcommand{\vvers}     {\mathop{\;{\longrightarrow}\;}}
\newcommand{\cvetroite} {\mathop{\;{\Longrightarrow}\;}}
\newcommand{\ieme}      {\hbox{i}^{\hbox{\smalltype\`eme}}}
\newcommand{\eqps}      {\, \buildrel \rm \hbox{\rm\smalltype p.s.} \over = \,}
\newcommand{\eqas}      {\,\buildrel\rm\hbox{\rm\smalltype a.s.} \over = \,}
\newcommand{\argmax}    {\hbox{{\rm Arg}}\max}
\newcommand{\argmin}    {\hbox{{\rm Arg}}\min}
\newcommand{\indep}{\perp\!\!\!\!\perp}
\newcommand{\abs}[1]{\left| #1 \right|}
\newcommand{\crochet}[2]{\langle #1 \,,\, #2 \rangle}
\newcommand{\espc}[3]   {E_{#1}\left(\left. #2 \right| #3 \right)}
\newcommand{\rang}{\hbox{rang}}
\newcommand{\rank}{\hbox{rank}}
\newcommand{\signe}{\hbox{signe}}
\newcommand{\sign}{\hbox{sign}}

\newcommand\hA{{\widehat A}}

% les lettres calligraphiques
%
\newcommand{\AAA}{{\cal A}}
\newcommand{\BB}{{\cal B}}
\newcommand{\CC}{{\cal C}}
\newcommand{\DD}{{\cal D}}
\newcommand{\EE}{{\cal E}}
\newcommand{\FF}{{\cal F}}
\newcommand{\GG}{{\cal G}}
\newcommand{\HH}{{\cal H}}
\newcommand{\II}{{\cal I}}
\newcommand{\JJ}{{\cal J}}
\newcommand{\KK}{{\cal K}}
\newcommand{\LL}{{\cal L}}
\newcommand{\NN}{{\cal N}}
\newcommand{\MM}{{\cal M}}
\newcommand{\OO}{{\cal O}}
\newcommand{\PP}{{\cal P}}
\newcommand{\QQ}{{\cal Q}}
\newcommand{\RR}{{\cal R}}
\newcommand{\sS}{{\cal S}}
\newcommand{\TT}{{\cal T}}
\newcommand{\UU}{{\cal U}}
\newcommand{\VV}{{\cal V}}
\newcommand{\WW}{{\cal W}}
\newcommand{\XX}{{\cal X}}
\newcommand{\YY}{{\cal Y}}

\renewcommand{\Box}{\diamond}
\newcommand{\tbullet}{$\bullet$}
\newcommand{\ot}{\leftarrow}
%\newcommand{\newblock}{}
%\newcommand{\carre}{\hfill$\Box$}
%\newcommand{\carreb}{\hfill\rule{0.25cm}{0.25cm}}
%
% utilitaires
%
%\newcommand{\dontforget}[1]
%{{\mbox{}\\\noindent\rule{1cm}{2mm}\hfill don't forget :
% #1 %\hfill\rule{1cm}{2mm}}\typeout{---------- don't forget : #1 ------------}}
%
%\newcommand{\note}[2]
%{ \noindent{\sf #1 \hfill \today}

%\noindent\mbox{}\hrulefill\mbox{}
%\begin{quote}\begin{quote}\sf #2\end{quote}\end{quote}
%\noindent\mbox{}\hrulefill\mbox{}
 
\vspace{1cm}
\newcommand{\be}{\begin{equation}}
\newcommand{\bea}{\begin{eqnarray}}
\newcommand{\eea}{\end{eqnarray}}

\newcommand{\beq}[1]{\begin{equation}\label{#1}}
\newcommand{\eeq}{\end{equation}}
\newcommand{\req}[1]{(\ref{#1})}
\newcommand{\beqn}[1]{\begin{eqnarray}\label{#1}}
\newcommand{\eeqn}{\end{eqnarray}}
\newcommand{\beaa}{\begin{eqnarray*}}
\newcommand{\eeaa}{\end{eqnarray*}}

\newcommand{\eq}[1]{(\ref{#1})}

\newcommand{\rond}[1]     {\stackrel{\circ}{#1}}
\newcommand{\rondf}       {\stackrel{\circ}{\FF}}
\newcommand{\point}[1]     {\stackrel{\cdot}{#1}}
\newcommand{\eps}{\varepsilon}
\newcommand\relatif{{\rm \rlap \tau\kern 3pt \tau}}
\newcommand{\limsupe}{\limsup_{\varepsilon \rightarrow 0}}
\newcommand{\liminfe}{\liminf_{\varepsilon \rightarrow 0}}
\newcommand{\lime}{\lim_{\varepsilon \rightarrow 0}}
\newcommand{\Aep}{A^\varepsilon}
\newcommand{\Bep}{B^\varepsilon}
\newcommand{\bX}{{\bar X}}
\newcommand{\tX}{{\tilde X}}
\newcommand{\hX}{{\hat X}}
%%%%%%%%%%%%%%%%%%%%%%%%%%%%%%%%%%%%%%%%%%%%%%%%%%%%%
\maketitle
\section*{Abstract}
Consider the standard, one dimensional, 
nonlinear filtering problem for diffusion processes observed in
small additive white noise:
$dX_t=b(X_t)dt+ dB_t\,,
dY_t^\eps=\gamma(X_t)dt+\eps dV_t\,,$ where $B_\cdot, V_\cdot$
are standard independent Brownian motions. Denote by 
$q^\eps_1(\cdot)$ the density of the
law of $\Xi_1$ conditioned on 
$\sigma(Y_t^\eps: 0\leq t\leq 1)$. We provide 
``quenched" large deviation estimates for 
the random family of measures $q^\eps_1(x)dx$: there exists a continuous, 
explicit
mapping $\bar \JJ~: \reels^2\to\reels$ such that for almost 
all $B_\cdot,V_\cdot$, $\bar \JJ(\cdot,X_1)$ is a good rate function and
for any measurable $G\subset \reels$, 
$$-\inf_{x\in G^o} \bar \JJ(x,X_1) \leq 
\liminfe \eps \log \int_G q_1^\eps(x) dx
\leq 
\limsupe \eps \log \int_G q_1^\eps(x) dx\nonumber
\leq 
-\inf_{x\in \bar G} \bar \JJ(x,X_1)\,.$$

\vspace{0.5cm}
\section{Introduction and statement of results}
Consider the following one dimensional filtering problem, where the signal process
$X_\cdot$ and the observation process $Y_\cdot^\eps$, parametrized by a
``small noise intensity'' $\eps$, are 
\bea
\label{intro-3}
dX_t&=& b(X_t)dt+ dB_t\,, 
\quad X_0\sim p_0(\cdot)\nonumber \\
dY_t^\eps&=&h(X_t)dt+\eps dV_t\,.
\eea
Here, $B_\cdot, V_\cdot$ are independent standard one dimensional
Brownian motions, and the functions $b,h,p_0$
satisfy the assumptions\footnote{Due to the one-dimensional nature of
our model, no generality is lost in assuming the diffusion coefficient 
of the signal process to be one. Indeed, if 
the signal process satisfies
$d\Xi_t=\beta(\Xi_t)dt+ \sigma(\Xi_t)dB_t,$ with $\sigma$ uniformly bounded 
away from zero,
then the ransformation 
$X_t=\bar G(\Xi_t)$, with
$\bar G(x)=\int_0^x (1/\sigma)(u) du\,,$
allows one to rewrite the problem in the form (\ref{intro-3}).
}
$$
\begin{array}{ll}
(A-1)& b,h,b',h'
 \,\mbox{are Lipschitz functions}\\
(A-2)& h'(\cdot)\geq h_0>0\\
(A-3)& |\log p_0(x)-\log p_0(y)|\leq c(1+|x|+|y|)|x-y|\,, x,y\in \reels\,,
\quad p_0 \, \mbox{\rm is uniformly bounded}\,.
\end{array}
$$
%Thus, in what follows, we consider only the filtering
%problem (\ref{intro-3}) satisfying (A-1)--(A-3). 
For technical reasons, we need to 
impose the following additional restriction:
$$
\begin{array}{ll}
(A-4)& h'b, h'h, h'', hb  \,\mbox{ are Lipschitz functions} \,, 
 \mbox{ and} \,
\lim_{|x|\to\infty}  h''(x)=0\,.
\end{array}
$$
$(A-4)$ implies that outside large compacts,  
the observation function $h$
function is essentially linear.
Let $\Omega_1=\Omega_2=C([0,1];\reels)$, 
$\Omega=\Omega_1\times \Omega_2$,
$\FF_i$ the Borel $\sigma$-algebra on $\Omega_i$, $i=1,2$,
$\FF$ the Borel $\sigma$-algebra on $\Omega$; let
$P_1,P_2$ denote the Wiener measure on $\Omega_1,\Omega_2$, and
$P=P_1\otimes P_2$. We define $B_t(\omega)=\omega_1(t)$, $V_t(\omega)=\omega_2(t)$, $0\le t\le1$. The pair $(B,V)$ is then distributed according to $P$. The 
solution $(X,Y^\eps)$ of
the  SDE (\ref{intro-3}) is then an $\FF$-measurable,
$C([0,1];\reels^2)$--valued, random variable.

Let $\mu_t^\eps(\cdot)$ denote the conditional law of $X_t$
conditioned on ${\cal Y}_t^\eps=\sigma\{Y^\eps_s,\, 0\le s\le t\}$, 
which we consider as an $\FF$-measurable map
from $\Omega$ to $M_1(\reels)$,
 the space of probability measures on $\reels$.
Note that $\mu_t^\eps$ is in fact measurable with respect to the
$\eps$-dependent $\sigma$-algebra $\YY_t^\eps\subset \FF$. 

It is known that
$\mu^\eps_t$ is absolutely continuous, with $\mu^\eps_t(dx)=q^\eps_t(x)dx$, and that
 as $\eps\to 0$, the conditional law $\mu_1^{\eps}(dx)=q^{\eps}_1(x)dx$ of 
$X_1$ given ${\cal Y}_1^\eps$ converges to the
 Dirac measure $\delta_{X_1}$ (all these facts can be found, e.g., in
\cite{picard}). In
particular, $X_1$ is measurable with respect to the limiting $\sigma$--algebra 
${\cal Y}_1^0$, since $h$ is one--to--one. It is known from the results of Picard 
\cite{picard} that the conditional law $\mu_1^{\eps}$ has a variance of order $\eps$, and
can be well approximated by a Gaussian law, which is given by an extended Kalman filter.

Our goal in this paper is to establish a large deviations result in the following sense.
Let $G$ be a measurable subset of $\reels$.
By the above remarks, we know that 
on the event
$\{X_1\not\in\overline{G}\}$, 
$\mu_1^{\eps}(G)\to 0$, $P$-almost surely. 
It turns out that
it goes to zero at exponential speed, i.e. roughly
like $\exp[-c_1(G)/\eps]$. What is the value of $c_1(G)=
-\lim\eps\log\mu_1^{\eps}(G)$
(if this limit exists), the ``rate function'', which tells us at which speed the quantity
$P(X_1\in G\,|\,{\cal Y}_1^\eps)$ 
goes to zero, whenever $X_1\not\in \overline{G}$? Clearly 
$c_1(G)$ must
depend on $X_1$ (at least intuitively through its distance to $\overline{G}$), and we 
shall see that
this is indeed the case. There is no surprise in the fact that $c_1(\cdot)$ is random,
since it tells us at which exponential speed the random measures $\mu_1^{\eps}$ converge
to the random measure $\mu_1^0=\delta_{X_1}$. Our results show
that it does not depend on anything else, in the sense that conditionally on
$\sigma(X_1)$, it is $P$-almost surely constant.

We call our result ``quenched'' (borrowing
that terminology from the theory of random media), meaning that the randomness of the
observation 
process 
is frozen. One could also discuss a ``semi-quenched'' large deviations statement by
computing the $P_1$-almost sure limit (if it exists) of
$$\eps\log \int \int_G q^\eps_1(x+X_1)dx dP_2\,,$$
while an ``annealed'' large deviations result would describe the
asymptotic behaviour of
$$\eps\log E\int_G q^\eps_1(x+X_1)dx.$$
Finally, one could also consider large deviations questions
at the level of the conditional measure itself, for example
questions concerning the rate of decay of probabilities of
the form $P(q_1^\eps(x) dx\in A)$, with $A$ a measurable subset
of the space of probability measures on $\reels$.
We hope to study all these elsewhere. 

 Let us now state our result. Define 
$$\bar \JJ(x,X_1)= \int_{X_1}^x (h(y)-h(X_1)) dy\,.$$
Our main result is the following theorem. For standard definitions 
concerning the LDP, see \cite{DZ}. For a set $G\subset\reels$, we
denote by $G^o$ its interior and by $\bar G$ its closure.
\begin{e-theorem}
\label{theo-1}
Assume (A-1)--(A-4). Then the family of (random) probability measures
$q_1^\eps(x)dx$
satisfies a {\it quenched} LDP (on the space $\reels$ equipped with
the standard euclidean norm) with continuous, good
rate function $\bar \JJ(\cdot,X_1)$. That is, 
for any measurable set $G\subset \reels$,
\bea
\label{main-theo1}
-\inf_{x\in G^o} \bar \JJ(x,X_1)& \leq &
\liminfe \eps \log \int_G q_1^\eps(x) dx\nonumber
\\
&\leq &
\limsupe \eps \log \int_G q_1^\eps(x) dx\nonumber
\\
&\leq &
-\inf_{x\in \bar G} \bar \JJ(x,X_1)\,, \quad P-a.s.
\eea
In fact, we have the estimate, valid for any fixed
compact set
$K_0\subset \reels$,
\be
\label{main-theo2}
\lime \sup_{x\in K_0} |\eps \log q_1^\eps(x)  +\bar \JJ(x,X_1)|
=0\,,\quad P-a.s.
\end{equation}
\end{e-theorem}
(It will be obvious 
from the proof that the fixed time $1$ can be replaced by any 
fixed time $t\in (0,\infty)$,
 that is the statement of  Theorem \ref{main-theo1}
remains true  with $q_t^\eps$ and $X_t$ replacing $q_1^\eps$ 
and $X_1$).

\noindent
{\bf Remarks} {\it 1.
In the particular case $h(x)=x$, Theorem \ref{theo-1} can be 
deduced from the results of \cite{Zei}.
\\
2. The reader could wonder why is the statement 
(\ref{main-theo1}) equivalent to the large deviations
principle on $\reels$ for $P$-almost $\omega$, 
since in (\ref{main-theo1}),
 the null
set on which the statement does not hold true may depend on $G$. Note
however that once the inequalities in
(\ref{main-theo1}) hold true for each interval
$G=(a,b)$ on a set of full measure $\Omega_{a,b}$,
set 
$$\Omega' =\cap_{a,b\in \rationels} \Omega_{a,b}\,,$$
and conclude that $P(\Omega')=1$ while (\ref{main-theo1})
holds true for all $\omega\in \Omega'$ and all open
intervals $G$ with rational endpoints. Since the latter
are a base for the topology on $\reels$, one concludes 
(see e.g. \cite[Theorem 4.1.11]{DZ}) that the full LDP holds for
each $\omega\in \Omega'$.}

We conclude this introduction with some comments about previous work and
possible applications and
extensions of our
result. 
Our motivation for the study
of the large deviations of the optimal filter is  their
need in a variety of applications such as tracking (see \cite{zakai})
or the study of the filter memory length (see \cite{atar}).
In the one dimensional 
linear observation case studied in \cite{Zei}, precise pointwise 
estimates  can be derived by comparison with the linear filtering problem,
whose (Gaussian) solution is known explicitly. In contrast, here,
the main tool used in the 
proof of Theorem \ref{theo-1}
is the representation, due to Picard \cite{picard},
of the density $q_1^\eps$ in 
terms of an auxiliary sub-optimal filter, and the availability of 
good estimates on the performance of this suboptimal filter.
These results are not available in the general multi-dimensional case.
When they are, e.g. in the setup discussed in
\cite{picard1}, we believe our analysis can be 
carried through. Hence, while our result is presently limited to one dimension,
we expect that its multidimensional extension to the case where the dimensions of the
state and observation coincide, and the observation function is one--to--one, could be 
deduced from the results of 
\cite{picard1}. Extension to the case where the dimension of the observation is smaller
than the dimension of the state 
(which is the most relevant one for applications)
would require completely new additional ideas, since the result would be of a completely
different nature (the limiting measure is no longer necessarily 
a Dirac measure, and even when it is, the convergence to the Dirac measure 
is at different speeds for different coordinates).

We finally note that Hijab \cite{hijab} has derived a (path) 
quenched large deviations for the
conditional density for systems in which  both the signal and
the observation noises are small. This is related, by 
a time change, to looking at short times (of order $\eps T$) of 
the filtering equations
\begin{eqnarray*}
dX_t^\eps&=& \frac{1}{\eps} \bar b(X_t^\eps) dt+ dB_t\,,\quad X_0^\eps=x\\
dY_t^\eps&=& h(X_t^\eps)dt +\eps dV_t\,.
\end{eqnarray*}
(Hijab's results are not stated in this way, but are equivalent to 
the description given here. Note that his setup
is more general than ours in that it applies to the multi-dimensional
setup and allows for general regular diffusion coefficients).  
Hijab's results are not
directly comparable with the LDP we derive here 
because of 
the different time interval on which they apply, and also because of 
the different type of conditioning
(his statement looks at the conditional density as a continuous functional of the
observation trajectory, and considers the LDP when this trajectory
is frozen.  It is thus not directly applicable as a quenched statement). 

\noindent
{\bf Convention:} Throughout the paper, when relevant, we made explicit on what 
parameters do constants depend, even if the actual value of the constant may change 
from line to line. When nothing explicit is mentioned, i.e. a generic
constant $C$ is used, it is assumed that it may depend
on the trajectories $\{X_\cdot\}$, $\{V_\cdot\}$, but
not on $\eps$.
For $\infty>t>0$, we use the notation
$||f||_t=\sup_{s\leq t}|f(s)|$,  with $||f||:=||f||_{1/\eps}$.
Finally, we use $\theta^t$ to denote the shift operator,
e.g. 
$\theta^t \tilde m(\cdot)=\tilde m(t+\cdot)$.

\section{Picard's formulation and a path integral}

The filtering problem we are going to analyze is (\ref{intro-3}), 
and the assumptions $(A-1)$--$(A-4)$ will be
assumed to hold throughout the paper. We also note that since nothing
is changed (in terms of the filtering problem) by adding
a constant to the observation function $h$, we may and will assume
throughout the paper that $h(0)=0$.

It is known from the results of Picard \cite{picard} that the 
conditional law $q^{\varepsilon}_1(x) dx$ has a small variance, and 
that there exist finite dimensional filters that provide good 
approximations of the unknown state. We shall now recall the formula 
derived by Picard \cite{picard} for $q_1^{\varepsilon}(x)$, 
which 
was used there
%in \cite{picard} 
to study approximate filters. 
It will be an essential tool for our large deviation results.

Define the approximate filter
$$ dM_t^\eps=b(M_t^\eps) dt +\frac{1}{\eps} (dY_t^\eps-h(M_t^\eps)dt)
\,,$$ with $M_0^\eps=0$,
% chosen arbitrarily and independently
%of $\eps$, 
and let
$\bar m_s=M_{1-s}^\eps$ and $\tilde m_s=\bar m_{\eps s}$, $s\in [0,1/\eps]$.

One of the main contributions of Picard in \cite[Proposition 4.2]{picard} 
was to express the 
conditional density $q_1^{\varepsilon}(x)$ in terms of 
the law of an auxiliary  process $\{\bX^x_{1-t}, 0\leq t\leq 1\}$, 
which fluctuates backward in time, starting at time $1$  
from the position $x$, around the trajectory of the 
approximate filter $M^\eps_\cdot$.
Performing a time change and a Girsanov transformation, Picard's result
can be rewritten as follows\footnote{For completeness, and 
since the computations involved are somewhat
lengthy, we present the derivation in an appendix at the end of the paper}. 
Define the process
$$ d\tilde Z_s^{\eps,x}=\left[-h(\tilde Z_s^{\eps,x})+ \tilde m_s h'(\tilde Z_s^{\eps,x})-\eps
b(\tilde Z_s^{\eps,x})\right]ds+\sqrt{\eps} d\tilde W_s\,,\quad \tilde Z_0^{\eps,x}=x\,,$$
with $\widetilde{W}$ a standard Brownian motion, independent of $B_\cdot,
V_\cdot$. Throughout, we let
$\E$ and $\P$ denote expectations and probabilities 
with respect to the law of the Brownian motion $\widetilde W_\cdot$.
Then
a version of the conditional
density of $X_1$ given ${\cal Y}_1^\eps$ is given by 
\be
\label{eq-picard1}
q_1^\eps
(x)=\frac{\rho^\eps_1(x)}{\int_\reels \rho^\eps_1(x) dx}\,,
\end{equation}
where
\be
\label{eq-picard2}
\rho^\eps_1(x):=
e^{-F(x,\tilde m_0)/\eps}
\E\left[ \exp\left(I_\eps(\tilde Z^{\eps,x}_{1/\eps}, 
0)+\int_0^{1/\eps} g_1(\tilde
Z_s^{\eps,x},\tilde m_s) ds+
\frac1{\eps} \int_0^{1/\eps} 
g_2(\tilde Z_s^{\eps,x},\tilde m_s) ds\right)\right]
\,,
\end{equation}
and 
\begin{eqnarray*}
F(z,m)&=& 
 \int_0^z (h(y)-h(m)) dy-m h(z)+h(m) z\,,\\
I_\eps(z,m)&=&\log p_0(z)+ \frac1{\eps}F(z,m)\,,
%\int_0^z (h(y)-h(m)) dy\,,
\\
g_1(z,m)&=&-mh'(z)b(z)+mh''(z)/2+h(z) b(z) -h'(z)/2-\eps b'(z)-h(z) b(m)\,,\\
g_2(z,m)&=&
h(z)h(m)-h^2(m)/2-m h(z) h'(z)+m^2 h'(z)^2/2\,.
\end{eqnarray*}
Note that the assumptions $(A-1)-(A-4)$ ensure that, for each given $m$,
 $g_1(\cdot,m),g_2(\cdot,m)$
are Lipschitz functions with Lipschitz constant uniformly bounded
for $m$ in compacts.

It is important to note that above, and throughout the paper, expressions of
the form $\E(\cdot)$ may still be random, due to their possible dependence
in $B_\cdot,V_\cdot$. Thus, any equality between such expressions is
to be understood in an a.s. sense. We will not explicitely mention this
in what follows.

Equipped with (\ref{eq-picard2}), one is tempted to apply standard tools of
large deviations theory, viz. the large deviations principle for $\tilde Z_\cdot^{\eps,x}$ and
Varadhan's Lemma, to the analysis of the exponential rate of decay of the $\P$ expectation
in (\ref{eq-picard2}). This temptation is quenched when one realizes that in fact,
the rate of growth of $\rho_1^\eps$ is exponential in $1/\eps^2$,
and it is only after normalization that one can hope to obtain the relevant
$1/\eps$ asymptotics. This fact, unfortunately, makes the analysis slightly more
subtle. In the next section, we present several lemmas, whose proof is deferred to Section
\ref{proofs-lemmas}, and 
show how to deduce Theorem \ref{theo-1} from these lemmas. 
Before closing this section, however,
we state the following easy a-priori estimates. 
Recall that according to our
convention, $||X||_1=\sup_{s\leq 1}|X_s|$:

\begin{e-lemme}
\label{lemma-m}
$||X||_1<\infty$, $P$-a.s., 
$$ |||\tilde m|||:=\limsup_{\eps \to 0} \sup_{t\in [0,1/\eps]}
|\tilde m_t|<\infty, \quad P-a.s.,$$
and for 
%each $\eta<1/2$ and any deterministic sequence
$T_\eps=\log(1/\eps)$, $|||\tilde m_X|||:=\sup_{s\in [0,T_\eps]}
|\tilde m_s-X_1|$,
\begin{equation}
\label{eq-revmor1}
 \limsup_{\eps \to 0} |||\tilde m_X|||
=0,\quad P-a.s..
\end{equation}
Further, 
%for each $\eta<1/2$
there exists a constant $C_{V,X}$ depending only
on $\{X_\cdot,V_\cdot\}$  such that
$$  \sup_{s\in [0,T_\eps]} 
|\tilde m_s-X_1|\leq C_{V,X}/\sqrt{T_\eps}\,,\quad P-a.s.$$
\end{e-lemme}

\noindent
{\bf Proof of Lemma \ref{lemma-m}:} The statement that
$||X||_1<\infty$ is part of the statement concerning 
existence of solutions to the SDE (\ref{intro-3}).
Next, we prove that
\begin{equation}
\label{eqq-1}
\limsup_{\eps\to 0}\sup_{t\leq 1} |M_t^\eps|<\infty\,.
\end{equation}
Indeed, fix  constants $C=C(||X||_1)$  and $\eps_0$
such that $h(y)-h(x)+\sup_{\eps\leq \eps_0} \eps b(x) <0$
for all $x\geq C$ and $|y|\leq ||X||_1$
(this is always possible because  $b, h$ are Lipschitz and $h'>h_0$).
Define the stopping times
$\tau_0=0,\theta_0=0$ and
$$\tau_i=\inf\{t>\theta_{i-1}: M_t^\eps= C\}\,,
\theta_i=\inf\{t>\tau_{i}: M_t^\eps= C+1\}\,.$$
By definition, $M_t^\eps\leq C+1$ for $t\in[\tau_i,\theta_i]$ while, for
$t\in [\theta_i,\tau_{i+1}]$ it holds that
for all $\eps<\eps_0$,
$$M_t^\eps=M_{\theta_i}^\eps+\int_{\theta_i}^t
[b(M_s^\eps)+\frac{1}{\eps}(h(X_s)-h(M_s^\eps))]
ds+ V_t-V_{\theta_i}\leq
C+1+2||V||_1\,.$$
We conclude that $\sup_{t\leq 1} M_t^\eps\leq C+1+2||V||_1<\infty$
for all $\eps<\eps_0$. A similar argument shows that
$\inf_{t\le1}M^\eps_t\ge-(C+1+2\|V\|_1)$.

To see the stated convergence of $\tilde m_s$ to $X_1$, recall
that $X_t$ and $V_t$ are almost surely H\"older($\eta$) continuous, for
all $\eta<1/2$.
Fix $t_0=1-2\eps T_\eps$, $t_1=1-\eps T_\eps$,
$\delta_\eps=1/\sqrt{T_\eps}$, and write
$Y_t=M_t^\eps-X_1$. With these notations,
$$Y_t=Y_{t_0}+\int_{t_0}^t
 \left[b(M_s^\eps)+\frac{h(X_s)-h(X_1)}{\eps}
\right]ds+
\frac1{\eps}\int_{t_0}^t(h(X_1)-h(M_s^\eps))ds+(V_t-V_{t_0})\,.$$
By the first part of the lemma, it holds that 
$|Y_{t_0}|\leq C$. We first show
that for some $\tau\in (t_0,t_1)$ it holds
that $|Y_\tau|\leq \delta_\eps$. Indeed, assume without
loss of generality that $Y_{t_0}>\delta_\eps$.
Then,
by the H\"{o}lder property of $X_\cdot$ and $V_\cdot$, it holds that
$$
\sup_{t\in (t_0,t_1)}
|V_{t}-V_{t_0}|\leq C(\eps T_\eps)^{\eta},\quad
\sup_{t\in (t_0,t_1)}
|X_{t}-X_{t_0}|\leq C(\eps T_\eps)^{\eta}\,.$$
%&\leq &
%C_2\left[1+2T_\eps \eps^{\eta'}\right]+
%\frac{h}{\eps}\int_{t_0}^t(X_1-M_s^\eps)ds\,.
%\end{eqnarray*}
Hence, if a $\tau$ as defined above does not exist,
then necessarily, using the Lipschitz continuity of
$h$,
$$-C\leq C_1\eps T_\eps(1+\frac{(\eps T_\eps)^{\eta}}{\eps})
-h_0\delta_\eps T_\eps+ C_1 (\eps T_\eps)^{\eta}\,,$$
which is clearly impossible unless $\eps\ge\eps_0$ for some
$\eps_0>0$. Now, for $\tau<t\le1$ we claim that
it is impossible that $Y_t>2\delta_\eps$.
Indeed, let $\theta'=\inf\{\tau<t\le1:  Y_t=2\delta_\eps\}$.
Repeating the argument above, we now obain that
if such a $\theta'$ exists, it must hold that
for some $\theta<2\eps T_\eps$,
$$\delta_\eps\leq C_1 \theta+C_1\frac{\theta^{\eta+1}}{\eps}+
C_1 \theta^{\eta}- \frac{h_0\delta_\eps \theta}{\eps}\,,$$
which again is impossible, unless $\eps\ge\eps'_0$, for some $\eps'_0>0$.
The case of $Y_t<-2\delta_\eps$ for some $t>t_0$ being handled
similarly, 
the conclusion follows.
$\Box$

\section{Auxilliary Lemmas and Proof of Theorem \ref{theo-1}}
Set
$J_\eps(x):=\rho_1^\eps(x) e^{F(x,\tilde m_0)/\eps}$ and
\begin{equation}
\bar L_\eps(x,t)=
\exp\left(\int_0^{t} \left(
g_1(\tilde
Z_s^{\eps,x},\tilde m_s) +
\frac1{\eps}
g_2(\tilde Z_s^{\eps,x},\tilde m_s) \right)
ds\right)
\end{equation}
and
\begin{equation}
L_\eps(x,t)=
\exp\left(I_\eps(\tilde Z^{\eps,x}_{t}, 0)\right) \bar L_\eps(x,t)\,.
\end{equation}
Although both $\bar L_\eps(x,t)$ and $L_\eps(x,t)$ depend on the path
$\tilde m_\cdot$, we omit this dependence when no confusion
occurs, while $L_\eps(x,t,m_{\cdot})$ will denote the quantity $L_\eps(x,t)$ with
$\tilde{m}_\cdot$ replaced by $m_\cdot$, and similarly for 
$\bar L_\eps$.

The following are the auxilliary lemmas alluded to above.
The proof of the first, Lemma \ref{lem-LDP},
is standard, combining 
large deviations estimates for solutions of SDE's 
(see e.g. \cite[Theorem 2.13, Pg. 91]{azencott})
with Varadhan's lemma (see e.g. \cite[Theorem 4.3.1, Pg. 137]{DZ}),
and is omitted. 
%\cite{MNS} and \cite{??} for similar statements.
%{\bf Should add the REF. and be more precise (number of the therorem in Freidlin-Wentzell ?)}
\begin{e-lemme}[Finite horizon LDP]
\label{lem-LDP}
Fix $T<\infty$ and a compact $K\subset\subset \reels$.
Define
$$ I_T(x,z):= \sup_{\phi\in H^1: \phi_0=x,\phi_T=z}
\int_0^T g_2(\phi_s,X_1)ds-\frac12
\int_0^T \left[\dot{\phi}_s+h(\phi_s)-X_1h'(\phi_s)\right]^2ds $$
Then, uniformly in $x,z\in K$,
$P.-a.s.$,
$$ \limsup_{\delta\to 0}\limsup_{\eps\to 0}
\left| \eps\log 
\E\left[ \bar L_\eps(x,T)
{\bf 1}_{\{|\tilde Z_T^{\eps,x}-z|<\delta\}}
\right]
-
I_T(x,z)\right|=0\,.$$
\end{e-lemme}
It is worthwhile noting the following simpler representation of $I_T(x,z)$:
\begin{equation}
\label{eq-IT}
I_T(x,z)=
\sup_{\phi\in H^1: \phi_0=x,\phi_T=z}
\left[ X_1 (h(z)-h(x))-h(X_1)(z-x) -\frac12 \int_0^T \left[\dot{\phi}_s
-(h(X_1)-h(\phi_s))\right]^2ds\right]\,.
\end{equation}
From this representation, the following is immediate:
\begin{equation}
\label{eq-equil}
I_T(X_1,X_1)= 0\,,
\end{equation}
and,
with $V_T(x):=I_T(x,X_1)$, it holds that
\begin{equation}
\label{VTlim}
V_T(x)\to_{T\to \infty} -X_1h(x)+h(X_1) x
%-\int_{X_1}^x (h(X_1)-h(y)) dy\,.
\end{equation}
This, and standard large deviations considerations, give
\begin{e-corollaire}
\label{cor-LDP}
Uniformly in $x,z\in K$,
$P-a.s.$,
\begin{eqnarray*}
&& \limsup_{T\to\infty}
\limsup_{\delta\to 0}\limsup_{\eps\to 0}
\\&&
\left| \eps\log 
\E\left[ \bar L_\eps(x,T)
{\bf 1}_{\{|\tilde Z_T^{\eps,x}-z|<\delta/2\}}
{\bf 1}_{\{|\tilde Z_{T/2}^{\eps,x}-X_1|<\delta/2\}}
\right]
-
h(X_1)x+h(x)X_1-I_{T/2}(X_1,z)\right|\\
&=& \limsup_{T\to\infty}
\limsup_{\delta\to 0}\limsup_{\eps\to 0}
\\&&
\left| \eps\log 
\E\left[ \bar L_\eps(x,T)
{\bf 1}_{\{|\tilde Z_T^{\eps,x}-z|<\delta/2\}}
{\bf 1}_{\{|\tilde Z_{T/2}^{\eps,x}-X_1|<\delta/2\}}
\right]
-
I_T(x,z)\right|\\
&=&  \limsup_{T\to\infty}
\limsup_{\delta\to 0}\limsup_{\eps\to 0}
\left| \eps\log 
\E\left[ \bar L_\eps(x,T)
{\bf 1}_{\{|\tilde Z_T^{\eps,x}-z|<\delta/2\}}
\right]
-
I_{T}(x,z)\right|
=0\,.
\end{eqnarray*}
\end{e-corollaire}

The key to the proof of Theorem \ref{theo-1} is a localization procedure that allows one
to restrict attention to compact (in time and space) subsets. A first coarse step in
that direction is provided by the next two lemmas. 
\begin{e-lemme}[Coarse localization 1]
\label{lem-coarse1}
For each $\eta>0$ there exists a constant \hfill\break
$M_1=M_1(|||\tilde m|||,\eta,|X_1|)$ and
$\eps_{00}=\eps_{00}(|||\tilde m|||,\eta,|X_1|)$ such that for all $\eps<\eps_{00}$,
\begin{equation}
\label{eq-xcut}
\int \rho_1^\eps(x){\bf 1}_{\{|x|>M_1/\sqrt{\eps}\}} dx 
\leq e^{-\eta/\eps} \inf_{|x|<1}  \rho_1^\eps(x)
\leq e^{-\eta/\eps} \int
 \rho_1^\eps(x){\bf 1}_{\{|x|\leq M_1/\sqrt{\eps}\}} dx 
\,, \quad P-a.s.
\end{equation}
\end{e-lemme}
\begin{e-lemme}[Coarse localization 2]
\label{lem-coarse2}
For each $\eta>0$ and $M_1$, $\eps_{00}$
 as in Lemma \ref{lem-coarse1},
there exist constants
$M_i=M_i(|||\tilde m|||,\eta,|X_1|)$,
$i=2,3$, with $M_3\leq M_2$, and
$\eps_0=\eps_0(|||\tilde m|||,\eta,|X_1|)<\eps_{00}$ such that for all $\eps<\eps_0$,
 uniformly in $|x|\leq M_1/\sqrt{\eps}$,
\begin{equation}
\label{eq-trajcut}
J_\eps(x)\leq 2 
\E\left[ L_\eps(x,1/\eps)
{\bf 1}_{\{||\tilde Z^{\eps,x}||\leq M_3/\eps\}}
\right]
\,,
\end{equation}
and uniformly in $|z|\leq M_3/\eps$, $T<1/\eps$,
\begin{equation}
\label{eq-kijung}
\E\left[ L_\eps(z,1/\eps-T,\theta^T\tilde m)\right]
\leq 2 \E\left[ L_\eps(z,1/\eps-T,\theta^T\tilde m)
{\bf 1}_{\{||\tilde Z^{\eps,z}||_{1/\eps-T}\leq M_2/\eps\}}\right]
\end{equation}
\end{e-lemme}
The following comparison lemma is also needed:
%uses the 
%coarse localization Lemma \ref{lem-coarsecomp}.
\begin{e-lemme}
\label{lem-comp}
There exists a function 
$g:\reels_+\mapsto\reels_+$, depending on $|||\tilde m|||, |X_1|,\eta$
only, with $g(\delta)\to_{\delta\to 0} 0$,
and an $\eps_1=\eps_1(|||\tilde m|||,X_1,\eta)<\eps_0$
such that
for all $\eps<\eps_1$, $t\in[1/2\eps,1/\eps]$,
 and 
$|x|, |y|\leq  M_3/ {\eps}$,
$|x-y|<\delta$,
\begin{equation}
\label{eq-comp1}
\eps \log\left(\frac{\E (L_\eps(x,t, \theta^{1/\eps-t}\tilde m)
{\bf 1}_{\{||\tilde Z^{\eps,x}||_t\leq M_2/\eps\}})
}{\E (L_\eps(y,t, \theta^{1/\eps-t}\tilde m)
{\bf 1}_{\{||\tilde Z^{\eps,y}||_t\leq M_2/\eps\}})
}\right) \leq g(\delta)\,,
\end{equation}
and there exists a constant $C_1(|||\tilde m|||,X_1,\eta)$ such that 
for all $\eps<\eps_1$,
\begin{equation}
\label{eq-comp2}
%\!\!\!\!\!\!\!\!\!
%\!\!\!\!\!\!\!\!\!
%\limsup_{\eps\to 0} 
%\sup_{|x|<M_2/\eps} 
\sup_{t\in [1/2\eps,1/\eps]}
\eps\left | \log \left( \frac{
\E\left[ L_\eps(x,t, \theta^{1/\eps-t}\tilde m)
{\bf 1}_{\{||\tilde Z^{\eps,x}||_t\leq M_2/\eps\}}
\right]}
{
\E\left[ L_\eps(X_1,t, \theta^{1/\eps-t}\tilde m)
{\bf 1}_{\{||\tilde Z^{\eps,X_1}||_t\leq M_2/\eps\}}
\right]}
\right)\right|
\leq C_1(1+ |x|)\,.
\end{equation}
\end{e-lemme}
The last step needed in order to carry out the localization procedure is the following
\begin{e-lemme}[Localization]
\label{lem-local}
Fix a sequence $T_\eps$ as in lemma \ref{lemma-m}. Then
there 
exist  constants $C_i=C_i(|||\tilde m|||,M_1,M_2,M_3,X_1)>0$, $i\geq 2$,
and a constant $\eps_2=\eps_2(|||\tilde m|||,M_1,M_2,M_3,X_1)<
\eps_1$,
such that for all $\eps<\eps_2$, 
$|x|\leq M_1/\sqrt{\eps}$, $|z|\leq M_3/\eps$,
$\delta<1$, and $1\leq T\leq T_\eps$,
\begin{equation} 
\label{eq-110}
\E\left[\bar L_\eps(x,T)
{\bf 1}_{\{|\tilde Z_T^{\eps,x}-z|<\delta\}}
{\bf 1}_{\{||\tilde Z^{\eps,x}||_T\leq M_3/\eps\}}
\right]
\leq \exp\left(\frac{C_2}{\eps}-\frac{C_3(|z|-|x|))_+^2}{\eps}
+\frac{C_4(|x|+|z|)}{\eps}\right)\,,
\end{equation}
and, uniformly for $|z-X_1|<1, |x-X_1|<1$,
\begin{equation} 
\label{eq-111}
\E\left[\bar L_\eps(x,T)
{\bf 1}_{\{|\tilde Z_T^{\eps,x}-z|<\delta\}}
{\bf 1}_{\{||\tilde Z^{\eps,x}||_T\leq M_3/\eps\}}
\right]
\geq 
\exp\left(-\frac{C_2}{\eps}\right)\,.
\end{equation}
\end{e-lemme}
We may now proceed to the proof of Theorem \ref{theo-1}, as
a consequence of the above Lemmata. 
Fix an $\eta>0$ as in Lemma
\ref{lem-coarse1}, and for $\delta>0$, $T>0$ to be chosen below,
with $T<T_\eps$, $T_\eps$ as in Lemma \ref{lemma-m}, 
 define
\begin{eqnarray}
\label{eq-revmor2}
\tilde J_\eps(x)&=&\E (L_\eps(x,1/\eps) 
{\bf 1}_{\{||\tilde Z^{\eps,x}||_T\leq M_3/\eps,
||\tilde Z^{\eps,x}||\leq M_2/\eps\}}) \nonumber \\
&=&
\sum_{i=-M_3/\eps \delta}^{M_3/\eps \delta}
\E (L_\eps(x,1/\eps) {\bf 1}_{\{||\tilde Z^{\eps,x}||\leq M_2/\eps, 
||\tilde Z^{\eps,x}||_T\leq M_3/\eps, 
|\tilde Z^{\eps,x}_T - i\delta|\leq \delta/2\}})\nonumber \\
&=:& 
\sum_{i=-M_3/\eps \delta}^{M_3/\eps \delta}
\tilde 
J_{\eps,T}(x,i\delta)\,.\end{eqnarray}
Set
${\cal Z}_T^{\eps,x}=\sigma(\tilde Z_t^{\eps,x}, t\leq T)$.
Using the Markov property, and the fact that $M_3<M_2$,
one may write, for $|z|<M_3/\eps$,
\begin{equation}
\label{eq-rev1}
\tilde J_{\eps,T}(x,z)=
\E\left[ 
\bar L_\eps(x,T)
{\bf 1}_{\{|\tilde Z_T^{\eps,x}-z|\leq
\delta/2\}}
{\bf 1}_{\{||\tilde Z^{\eps,x}||_T\leq M_3/\eps\}}
\E\left( L_\eps(\tilde
Z_T^{\eps,x},1/\eps-T, \theta^T \tilde m)
{\bf 1}_{\{||\tilde Z^{\eps,x}||\leq M_2/\eps\}}
\,|\,{\cal Z}_T^{\eps,x}\right)
\right]\,.
\end{equation}
Applying (\ref{eq-comp1}) and the Markov property, it follows that
on the event $\{|\tilde Z^{\eps,x}_T-z|\leq \delta/2\}\cap
\{||\tilde Z^{\eps,x}||_T\leq M_3/\eps\}$, 
%and denoting
%by ${\cal Z}_T^{\eps,x}=\sigma(\tilde Z_t^{\eps,x}, t\leq T)$, 
one has
for $\eps<\eps_1$, and $|x|\leq M_1/\sqrt{\eps}$,
$|z|\leq M_3/\eps$,
\begin{eqnarray*}
\E\left( L_\eps(\tilde
Z_T^{\eps,x},1/\eps-T, \theta^T \tilde m)
{\bf 1}_{\{||\tilde Z^{\eps,x}||\leq M_2/\eps\}}
\,|\,{\cal Z}_T^{\eps,x}\right)
&=&
\E\left( L_\eps(\tilde
Z_T^{\eps,x},1/\eps-T, \theta^T \tilde m)
{\bf 1}_{\{\sup_{T\leq t\leq 1/\eps}|\tilde Z^{\eps,x}_t|\leq M_2/\eps\}}
\,|\,{\cal Z}_T^{\eps,x}\right)\\
&
\leq&
e^{g(\delta)/\eps}
\E\left( L_\eps(z,
1/\eps-T, \theta^T \tilde m)
{\bf 1}_{\{\sup_{0\leq 1/\eps-T}|\tilde Z^{\eps,z}_t|\leq M_2/\eps\}}
\right)\\
&=&
e^{g(\delta)/\eps}
\E\left( L_\eps(z,
1/\eps-T, \theta^T \tilde m)
{\bf 1}_{\{||\tilde Z^{\eps,z}||_{1/\eps-T}\leq M_2/\eps\}}
\right)\,.
\end{eqnarray*}
Substituting in (\ref{eq-rev1}), one concludes that  for all
$\eps<\eps_1$, and $|x|\leq M_1/\sqrt{\eps}$,
$|z|\leq M_3/\eps$,
\begin{eqnarray}
\label{turkey1}
\tilde J_{\eps,T}(x,z)e^{-g(\delta)/\eps}
&\leq& 
\E\left[\bar L_\eps(x,T)
{\bf 1}_{\{|\tilde Z_T^{\eps,x}-z|
\leq \delta/2\}}
{\bf 1}_{\{||\tilde Z^{\eps,x}||_T\leq M_3/\eps\}}\right] \\
&&\quad 
\cdot \E \left[L_\eps(z,1/\eps-T,\theta^T\tilde m)
{\bf 1}_{\{||\tilde Z^{\eps,z}||_{1/\eps-T}\leq M_2/\eps\}}\right]
:= \hat J_{\eps,T}(x,z)
\leq
\tilde J_{\eps,T}(x,z)e^{g(\delta)/\eps}
\,.
\nonumber
\end{eqnarray}
Next, using (\ref{eq-comp2}) in the first inequality and Lemma
\ref{lem-local} in the second, it follows that for all $\eps<
\eps_2$, and $T\in (1,T_\eps)$, $T_\eps$ as in Lemma
\ref{lemma-m}, and some constants $C_i$ independent of $T$,$\eps$,
\begin{eqnarray}
\label{eq-rev10}
 \hat J_{\eps,T}(x,z)&\leq &
\E\left[
\bar L_\eps(x,T)
{\bf 1}_{\{|\tilde Z_T^{\eps,x}-z|
\leq \delta/2\}}
{\bf 1}_{\{||\tilde Z^{\eps,x}||_T\leq M_3/\eps\}}
\right]\nonumber\\
&& \quad\cdot
\E \left[L_\eps(X_1,1/\eps-T,\theta^T \tilde m)
{\bf 1}_{\{||\tilde Z^{\eps,X_1}||_{1/\eps-T}\leq M_2/\eps\}}\right]
e^{C_1(|z|+1)/\eps}\nonumber\\
&\leq &
\exp\left(\frac{C_2}{\eps}-\frac{C_3(|z|-|x|)_+^2}{\eps}
+\frac{C_5(|x|+|z|)}{\eps}\right)\nonumber\\
&&\quad \cdot
\E \left[L_\eps(X_1,1/\eps-T,\theta^T \tilde m)
{\bf 1}_{\{||\tilde Z^{\eps,X_1}||_{1/\eps-T}\leq M_2/\eps\}}\right]\,.
\end{eqnarray}
Similarly, for all $\eps<\eps_2$,  and $|x-X_1|\leq 1$,
$|z-X_1|\leq 1$,
\begin{equation}
\label{eq-rev11}
 \hat J_{\eps,T}(x,z)
\geq \exp\left(-\frac{C_2}{\eps}\right)
\E \left[L_\eps(X_1,1/\eps-T,\theta^T \tilde m)
{\bf 1}_{\{||\tilde Z^{\eps,X_1}||_{1/\eps-T}
\leq M_2/\eps\}}\right]\,.
\end{equation}
We next note that due to the 
quadratic growth of $F(x,X_1)$ as $|x|\to \infty$,
there exists a compact set ${\cal K}_1$, depending
on $|||m|||,X_1,\eta, C_i$ only,
such that 
\begin{equation}
\label{eq-rev12}
\sup_{(x,z)\in ({\cal K}_1\times {\cal K}_1)^c}
\frac{C_2}{\eps}-\frac{C_3(|z|-|x|)_+^2}{\eps}
+\frac{C_5(|x|+|z|)}{\eps} -\frac{F(x,X_1)}{\eps}
\leq -\frac{F(X_1,X_1)}{\eps} -\frac{C_2}{\eps}\,.
\end{equation}
Thus, using (\ref{eq-rev10}) in the first inequality, 
(\ref{eq-rev12}) in the second, and
(\ref{eq-rev11}) in the third,
\begin{eqnarray}
\label{eq-rev13}
&&\sup_{|x|\leq M_1/\sqrt{\eps},|z|\leq M_3/\eps,
(x,z)\in ({\cal K}_1\times {\cal K}_1)^c}
\hat J_{\eps,T}(x,z) e^{-F(x,X_1)/\eps}
\nonumber \\
&\leq &
\E \left[L_\eps(X_1,1/\eps-T,\theta^T \tilde m)
{\bf 1}_{\{||\tilde Z^{\eps,X_1}||_{1/\eps-T}
\leq M_2/\eps\}}\right]
\nonumber\\
&&\quad \cdot 
\sup_{|x|\leq M_1/\sqrt{\eps},|z|\leq M_3/\eps,
(x,z)\in ({\cal K}_1\times {\cal K}_1)^c}
\exp\left(\frac{C_2}{\eps}-\frac{C_3(|z|-|x|)_+^2}{\eps}
+\frac{C_5(|x|+|z|)}{\eps}-\frac{F(x,X_1)}{\eps}\right)\nonumber\\
&\leq&
\E \left[L_\eps(X_1,1/\eps-T,\theta^T \tilde m)
{\bf 1}_{\{||\tilde Z^{\eps,X_1}||_{1/\eps-T}
\leq M_2/\eps\}}\right]
%\nonumber\\
%&&\quad \cdot 
\exp\left(-\frac{C_2}{\eps}-\frac{F(X_1,X_1)}{\eps}
\right)\nonumber\\
&\leq & \hat J_{\eps,T}(X_1,X_1) e^{-F(X_1,X_1)/\eps}\,.
\end{eqnarray}
It follows by substituting (\ref{eq-rev13})
into (\ref{turkey1}) that
for all $\eps$ small enough,
and any $T\in (0,T_\eps)$,
\begin{equation}
\label{eq-compact}
\sup_{|x|\leq M_1/\sqrt{\eps},|z|\leq M_3/\eps}
\tilde J_{\eps,T}(x,z) e^{-F(x,X_1)/\eps}
\leq e^{2 g(\delta)/\eps}
\sup_{x\in {\cal K}_1, z\in {\cal K}_1}
\tilde J_{\eps,T}(x,z) e^{-F(x,X_1)/\eps}
\,.
\end{equation}
We may, by enlarging ${\cal K}_1$ if necessary, assume also that
%\begin{equation} \inf_{x\not\in {\cal K}_1} \bar{\cal I}(x,X_1)
%\geq 4\eta\,.
%\end{equation}
$[-1,1]\subset {\cal K}_1$.
With $\eta$ and ${\cal K}_1$,
as above, choose next $T$ large enough, $\delta$ small
enough (with $g(\delta)<\eta/8$)
and $\eps_3(\delta,T,\eta,|||\tilde m|||,|||\tilde m_X|||,X_1)$
$<\eps_2$
 such
that, 
for all $\eps<\eps_3$:
\begin{itemize}
\item 
The errors in the expression in
Corollary \ref{cor-LDP} and in (\ref{VTlim}) are each
bounded above by $\eta/8$,
uniformly in
$x, z\in {\cal K}_1$.
\item $|F(x,\tilde m_0)-F(x,X_1)|\leq \frac{\eta}{8}$, uniformly
in
$x\in {\cal K}_1$ (which is possible by Lemma \ref{lemma-m} and the 
uniform continuity of
$F(x,\cdot)$ for $x$ in compacts).
\item $\eps \log 2 \leq \frac{\eta}{8}$.
\item $\eps\log (2M_3/\eps\delta)\leq \frac{\eta}{8}$.
\end{itemize}
Hence,
%Using the Markov property again, Lemma \ref{lemma-m} and Corollary
%\ref{cor-LDP},
%one then gets that 
for 
$x\in {\cal K}_1$,
and 
all $\eps<\eps_3$,
\begin{eqnarray}
\label{eq-fin1}
\eps\log \rho_1^\eps(x)&=&
-F(x,\tilde m_0) 
+\eps\log \E(L_\eps(x,1/\eps))\quad \quad \quad{\mbox{\tt by  (\ref{eq-picard2})}}
\nonumber \\
&\leq & 
-F(x,\tilde m_0) +
\eps\log \E(L_\eps(x,1/\eps){\bf 1}_{\{||\tilde Z^{\eps,x}||\leq M_3/\eps
\}})+\eps \log 2\quad \mbox{\tt by  (\ref{eq-trajcut})}
\nonumber \\
&\leq &-F(x,X_1)+
\eps\log \E(L_\eps(x,1/\eps){\bf 1}_{\{||\tilde Z^{\eps,x}||\leq M_3/\eps
\}})+\frac{\eta}{4}\quad \mbox{\tt by  $\eps<\eps_3$}\nonumber \\
&\leq &-F(x,X_1)+
\eps\log \tilde J_\eps(x)
+\frac{\eta}{4}\quad \mbox{\tt by  (\ref{eq-revmor2})}%\,.
\nonumber \\
%\end{eqnarray}
%Therefore,
%for 
%$x\in {\cal K}_1$,
%and 
%all $\eps<\eps_3$,
%\begin{eqnarray}
%\label{eq-fin1}
%\eps\log \rho_1^\eps(x)
&\leq &-F(x,X_1)+
\eps\log \sup_{z\in {\cal K}_1}
\tilde J_{\eps,T}(x,z)
+\frac{\eta}{2}\quad \mbox{\tt by  (\ref{eq-revmor2}) and
(\ref{eq-compact})}\nonumber \\
&\leq &-F(x,X_1) 
+\sup_{z\in {\cal K}_1}
\left[\eps \log \E(\bar L_\eps(x,T)
{\bf 1}_{\{|\tilde Z_{T}^{\eps,x}-z|
\leq \delta/2\}}
%{\bf 1}_{\{||\tilde Z^{\eps,x}||_T\leq M_3/\eps\}})
)
\right.\nonumber\\
&&\quad \left.+\eps\log \E(L_\eps(z,1/\eps-T,\theta^T \tilde m)
{\bf 1}_{\{||\tilde Z^{\eps,z}||_{1/\eps-T}\leq M_2/\eps\}})\right]
+\frac{5\eta}{8}\quad\mbox{\tt by (\ref{turkey1})}\nonumber \\
&\leq &-F(x,X_1) 
+\sup_{z\in {\cal K}_1}
\left[h(X_1)x-h(x)X_1+I_{T/2}(X_1,z)
%\eps \log \E(\bar L_\eps(x,T)
%{\bf 1}_{\{|\tilde Z_{T/2}^{\eps,x}-X_1|
%\leq \delta/2\}}
%{\bf 1}_{\{|\tilde Z_{T}^{\eps,x}-z|
%\leq \delta/2\}}
%)
\right.\nonumber\\
&&\quad \left.+\eps\log \E(L_\eps(z,1/\eps-T,\theta^T \tilde m)
{\bf 1}_{\{||\tilde Z^{\eps,z}||_{1/\eps-T}\leq M_2/\eps\}})\right]
+\frac{7\eta}{8}\quad\mbox{\tt by Corollary(\ref{cor-LDP})}\nonumber \\
%&\leq &\eps\log \E(\bar L_\eps(x,T/2)
%{\bf 1}_{\{|\tilde Z_{T/2}^{\eps,x}-X_1|
%\leq \delta/2\}})\nonumber \\
%&& \;\;
%+
%\sup_{z\in {\cal K}_1}
%\left[\eps\log \E(\bar L_\eps(X_1,T/2,\theta^{T/2} \tilde m)
%{\bf 1}_{\{|\tilde Z_{T/2}^{\eps,X_1}-z|
%\leq \delta/2\}})+
%\eps\log \E(L_\eps(z,1/\eps-T,\theta^T \tilde m))\right]\nonumber\\
&\leq &
h(X_1)x-h(x)X_1-F(x,X_1)+\eta +
%\\
%&& \;\;\;
\sup_{z\in {\cal K}_1}
\left[I_{T/2}(X_1,z)+
%\eps\log \E(\bar L_\eps(X_1,T/2,\theta^{T/2}\tilde m)
%{\bf 1}_{\{|\tilde Z_{T/2}^{\eps,X_1}-z|
%\leq \delta/2\}})+
\eps\log \E(L_\eps(z,1/\eps-T,\theta^T \tilde m))\right]\nonumber 
\\
&=:&-\bar {\cal J}(x,X_1)+ \eta+ {\cal C}_\eps\,,
%\nonumber
\end{eqnarray}
where ${\cal C}_\eps$ depends only on $\eps$, and not on $x$, and
is defined by the last equality.
Similarly, for all $x\in {\cal K}_1$ and all $\eps<\eps_3$,
\begin{eqnarray}
\label{eq-fin2}
\eps\log \rho_1^\eps(x)&=&
-F(x,\tilde m_0) 
+\eps\log \E(L_\eps(x,1/\eps))\quad \quad \quad{\mbox{\tt by  (\ref{eq-picard2})}}
\nonumber \\
&\geq & 
-F(x,\tilde m_0) +
\eps\log \E(L_\eps(x,1/\eps){\bf 1}_{\{
||\tilde Z^{\eps,x}||_T\leq M_3/\eps,
||\tilde Z^{\eps,x}||\leq M_2/\eps
\}})
\nonumber \\
&\geq &-F(x,X_1)+
\eps\log \E(L_\eps(x,1/\eps){\bf 1}_{\{
||\tilde Z^{\eps,x}||_T\leq M_3/\eps,
||\tilde Z^{\eps,x}||\leq M_2/\eps
\}})-\frac{\eta}{4}\quad \mbox{\tt by  $\eps<\eps_3$}\nonumber \\
&= &-F(x,X_1)+
\eps\log \tilde J_\eps(x)
-\frac{\eta}{4}\quad \mbox{\tt by  definition}%\,.
\nonumber \\
&\geq &-F(x,X_1)+
\eps\log \sup_{z\in {\cal K}_1}
\tilde J_{\eps,T}(x,z)
-\frac{\eta}{4}\quad \mbox{\tt by  definition
}\nonumber \\
&\geq &-F(x,X_1) 
+\sup_{z\in {\cal K}_1}
\left[\eps \log \E(\bar L_\eps(x,T)
{\bf 1}_{\{|\tilde Z_{T}^{\eps,x}-z|
\leq \delta/2\}}
)
\right.\nonumber\\
&&\quad \left.+\eps\log \E(L_\eps(z,1/\eps-T,\theta^T \tilde m)
{\bf 1}_{\{||\tilde Z^{\eps,z}||_{1/\eps-T}\leq M_2/\eps\}})\right]
-\frac{5\eta}{8}\quad\mbox{\tt by (\ref{turkey1})}\nonumber \\
&\geq &-F(x,X_1) 
+\sup_{z\in {\cal K}_1}
\left[h(X_1)x-h(x)X_1+I_{T/2}(X_1,z)
\right.\nonumber\\
&&\quad \left.+\eps\log \E(L_\eps(z,1/\eps-T,\theta^T \tilde m)
{\bf 1}_{\{||\tilde Z^{\eps,z}||_{1/\eps-T}\leq M_2/\eps\}})\right]
-\frac{7\eta}{8}\quad\mbox{\tt by Corollary(\ref{cor-LDP})}\nonumber \\
&\geq &
h(X_1)x-h(x)X_1-F(x,X_1)-\eta +
\sup_{z\in {\cal K}_1}
\left[I_{T/2}(X_1,z)+
\eps\log \E(L_\eps(z,1/\eps-T,\theta^T \tilde m))\right]\nonumber 
\\
&=&-\bar {\cal J}(x,X_1)- \eta+ {\cal C}_\eps\,,
\end{eqnarray}
where ${\cal C}_\eps$ is the same as in (\ref{eq-fin1}).
%and hence
%\begin{equation}
%\label{eq-fin2}
%\eps\log \rho_1^\eps(x)\geq
%- \bar {\cal J}(x,X_1)- \eta+ {\cal C}_\eps\,.
%\end{equation}
Since
$\bar {\cal J}(\cdot,X_1)$ is continuous and 
$\bar {\cal J}(X_1,X_1)=0$, it follows from (\ref{eq-fin2}) that  
\begin{equation}
\label{eq-reveve2}
\liminf_{\eps\to 0}
\eps\log \int_\reels \rho_1^\eps(x) dx - {\cal C}_\eps
\geq  -2\eta\,.\end{equation}
On the other hand,  for $\eps<\eps_3$,
\begin{eqnarray}
\label{eq-reveve1}
\eps\log\int_\reels \rho_1^\eps(x)dx&\leq
&\eps \log(1+e^{-\eta/\eps})+
\eps\log\int_{|x|\leq M_1/\sqrt{\eps}}
\rho_1^\eps(x)dx \quad \mbox{\tt by Lemma \ref{lem-coarse1}}
\nonumber \\
&\leq&
\eps \log(1+e^{-\eta/\eps})+\eps\log 2+\eps\log\left(
\frac{2M_3}{\eps\delta}\right)
\nonumber \\
&&\quad +\sup_{|x|\leq M_1/\sqrt{\eps},|z|\leq M_3/\eps}
\eps\log \left(\tilde J_{\eps,T}(x,z)e^{-F(x,X_1)/\eps}\right)\quad
\mbox{\tt by Lemma \ref{lem-coarse2} and (\ref{eq-revmor2})}
\nonumber \\
&\leq &
\frac{5\eta}{8}
+\sup_{x,z\in {\cal K}_1}
\eps\log \left(\tilde J_{\eps,T}(x,z)e^{-F(x,X_1)/\eps}\right)\quad
\mbox{\tt by (\ref{eq-compact})}\nonumber\\
&\leq &
\frac{5\eta}{8}
+\eps\log\left(\sup_{x\in {\cal K}_1}\rho_1^\eps(x)\right)\nonumber\\
&\leq & 2\eta+{\cal C}_\eps-\inf_x \bar{\cal J}(x,X_1)=
2\eta+{\cal C}_\eps\quad\mbox{\tt by (\ref{eq-fin1}) and 
$\bar {\cal J}(x,X_1)\geq 0$}.
\end{eqnarray}
%using that $\bar {\cal J}(x,X_1)\geq 0$,
% with
%$\bar {\cal J}(X_1,X_1)=0$,  
%(\ref{eq-compact}) and Lemma \ref{lem-coarse1},
%it follows that 
%$$\lim_{\eps\to 0}\left[
%\eps\log \int_\reels \rho_1^\eps(x) dx
%- {\cal C}_\eps\right]=0\,.$$
Consider now an open ball $B(x_0,r)\subset \reels$. Then,
using (\ref{eq-reveve1}) in the first inequality,
and (\ref{eq-fin2}) in the last,
\begin{eqnarray*}
%\label{eq-lbld}
\liminf_{\eps\to 0} \eps\log \int_{B(x_0,r)}
q_1^\eps(x) dx&=&
\liminf_{\eps\to 0} [\eps\log \int_{B(x_0,r)}
\rho_1^\eps(x) dx
-\eps\log \int_\reels \rho_1^\eps(x)dx]\nonumber
\\ &\geq&
\liminf_{\eps\to 0} [\eps\log \int_{B(x_0,r)}
\rho_1^\eps(x) dx-{\cal C}_\eps-2\eta]\nonumber
\\ &\geq&
-
\bar{\cal J}(x_0,X_1)-3\eta\,.
\end{eqnarray*}
$\eta$ being arbitrary, one deduces that
\begin{equation}
\label{eq-lbld}
\liminf_{\eps\to 0} \eps\log \int_{B(x_0,r)}
q_1^\eps(x) dx\geq 
-
\bar{\cal J}(x_0,X_1)\,.
\end{equation}

To see the complementary upper bound for the ball
$B(x_0,r)$, enlarge ${\cal K}_1$ if necessary so that
$B(x_0,r)\subset {\cal K}_1$ (decreasing $\eps_3$ above
as a by product). Then, 
using (\ref{eq-reveve2}) in the first inequality,
and (\ref{eq-fin1}) in the last,
\begin{eqnarray*}
%\label{eq-lbld}
\limsup_{\eps\to 0} \eps\log \int_{B(x_0,r)}
q_1^\eps(x) dx&=&
\limsup_{\eps\to 0} [\eps\log \int_{B(x_0,r)}
\rho_1^\eps(x) dx
-\eps\log \int_\reels \rho_1^\eps(x)dx]\nonumber
\\ &\leq&
\limsup_{\eps\to 0} [\eps\log \int_{B(x_0,r)}
\rho_1^\eps(x) dx-{\cal C}_\eps+2\eta]\nonumber
\\ &\leq&
-
\sup_{x\in B(x_0,r)}
\bar{\cal J}(x,X_1)+3\eta+\limsup_{\eps\to 0}\eps\log (2r)\,.
\end{eqnarray*}
$\eta$ being arbitrary, the above, 
(\ref{eq-lbld}), and
the continuity of $\bar {\cal J}(\cdot,X_1)$ imply
that
$$
\lim_{r\to0}\limsup_{\eps\to 0}
\eps\log \int_{B(x_0,r)}
q_1^\eps(x) dx =
\lim_{r\to0}\liminf_{\eps\to 0}
\eps\log \int_{B(x_0,r)}
q_1^\eps(x) dx = \bar{\cal J}(x_0,X_1)\,.
$$
Next, \cite[Theorem 4.1.11]{DZ}, the above, Remark 2 following
Theorem \ref{theo-1}, and
the continuity of $\bar 
{\cal J}(\cdot,X_1)$ imply that
the weak LDP holds for the sequence of (random) measures
$\mu_1^\eps(dx)=q_1^\eps(x) dx$ on $\reels$. To prove the full large 
deviations principle, it remains, by
\cite[Lemma 1.2.8]{DZ}, to prove the exponential
tightness of the sequence $\mu_1^\eps$. That is, for each given
$L$ we must find a constant $C_L$ such that
\begin{equation}
\label{exp-tight}
\limsup_{\eps\to 0}
\eps
\log \int_{[-L,L]^c} q_1^\eps(x) ds <-L\,.
\end{equation}
Since the proof of (\ref{exp-tight}) uses some estimates
from the proof of Lemma \ref{lem-coarse1}, to avoid 
repetitions we postpone it to the end of Section \ref{proofs-lemmas}.

Finally, we note that (\ref{main-theo2}) is an immediate 
consequence of the estimates
(\ref{eq-fin1}), (\ref{eq-fin2}), 
(\ref{eq-reveve1}) and (\ref{eq-reveve2}).
$\Box$

\section{Proofs of auxilliary lemmas}
\label{proofs-lemmas}

Throughout this section, $C$ denotes a positive constant
that depends on $|||\tilde m|||, 
|||\tilde m_X|||, C_{V,X}, X$ only, and whose value may change 
from line to line.

\noindent
{\it Proof of Lemma \ref{lem-coarse1}}
The right inequality is a trivial consequence of the left one.
To prove the latter,
we first need an upper bound for the left hand side of
(\ref{eq-xcut}). 
A subsequent, easily derived 
 lower bound on the middle term  will conclude 
the proof.
Define the function
\begin{equation}
\label{minnts} H(x)=\int_0^x{h}(y)dy\,. 
%\text{ with }\overline{h}(x):=h(x)-h(0).
\end{equation}
We note that
\begin{equation*}
I_\eps(\tilde{Z}^{\eps,x}_{1/\eps},0)-\frac{F(x,\tilde{m}_0)}{\eps}
=\frac{1}{\eps}\left(H(\tilde{Z}^{\eps,x}_{1/\eps})-H(x)\right)
+\log p_0(\tilde{Z}^{\eps,x}_{1/\eps})+\frac{\tilde{m}_0h(x)
%h(\tilde{m}_0)x}{\eps}.
}{\eps}.
\end{equation*}
We first rewrite the $\tilde{Z}^{\eps,x}_t$ equation as
\[\tilde{Z}^{\eps,x}_t=x+\int_0^t[-h(\tilde{Z}^{\eps,x}_s)+
g(s,\tilde{Z}^{\eps,x}_s)]ds+\sqrt{\eps}
\tilde{W}_t,\]
and next deduce from It\^o's formula that
\begin{equation*}
 H(\tilde{Z}^{\eps,x}_{1/\eps})-H(x)=\int_0^{1/\eps}[-
h^2(\tilde{Z}^{\eps,x}_s)
 +(h g)(s,\tilde{Z}^{\eps,x}_s)+
\frac{\eps}{2} h'(\tilde{Z}^{\eps,x}_s)]ds
 +\sqrt{\eps}\int_0^{1/\eps} h(\tilde{Z}^{\eps,x}_s)d\tilde{W}_s.
\end{equation*}
It now follows from \eqref{eq-picard2} 
and the (uniform in $m$ in compacts)
linear growth of $g_1(z,m)$ and $g_2(z,m)$
in $z$, that for some $C$ (depending
on $|||\tilde m|||$ and $X$ only)
%(depending on 
%$\tilde m, X$) 
and all $\eps\le1$, $\delta>0$,
\begin{equation*}
\begin{split}
\rho_1^\eps(x)&\le\exp\left[ \frac{C}{\eps^2}+\frac{\tilde{m}_0h(x)
}{\eps} \right]
\left(\E\left[p_0(\tilde{Z}^{\eps,x}_{1/\eps})
\right]^{\frac{1+\delta}{\delta}}
\right)^{\frac{\delta}{1+\delta}}\\
&\quad\times
\left(\E\exp\left[\frac{1+\delta}{\sqrt{\eps}}
\int_0^{1/\eps}h(\tilde{Z}^{\eps,x}_s)d\tilde{W}_s
-\frac{1+\delta}{\eps}\int_0^{1/\eps}h^2(\tilde{Z}^{\eps,x}_s)ds+
\frac{C}{\eps}\int_0^{1/\eps}|\tilde{Z}^{\eps,x}_s|ds
\right]\right)^{\frac{1}{1+\delta}}.
\end{split}
\end{equation*}
Now provided $\delta<1$, $1+\delta>\frac{(1+\delta)^2}{2}$, 
and thus
there exists a $p>1$ and a $p'>0$ such that
$$ 1+\delta=\frac{p (1+\delta)^2}{2}+p'\,.$$
Thus, with $q=p/(p-1)$,
\begin{equation*}
\begin{split}
&\left(
\E\exp\left[\frac{1+\delta}{\sqrt{\eps}}
\int_0^{1/\eps}h(\tilde{Z}^{\eps,x}_s)d\tilde{W}_s
-\frac{(1+\delta)^2p}{2\eps}\int_0^{1/\eps}
h^2(\tilde{Z}^{\eps,x}_s)ds
-\frac{p'}{\eps}\int_0^{1/\eps}
h^2(\tilde{Z}^{\eps,x}_s)ds
+
\frac{C}{\eps}\int_0^{1/\eps}|\tilde{Z}^{\eps,x}_s|ds
\right]\right)^{\frac{1}{1+\delta}}\\
&\leq
\left(\E\exp\left[\frac{p(1+\delta)}
{\sqrt{\eps}}\int_0^{1/\eps}h(\tilde{Z}^{\eps,x}_s)d\tilde{W}_s
-\frac{(1+\delta)^2p^2}{2\eps}\int_0^{1/\eps}
h^2(\tilde{Z}^{\eps,x}_s)ds
\right]\right)^{\frac{1}{p(1+\delta)}}\\
&\quad \times
\left(\E\exp\left[
-\frac{p'q}{2\eps}\int_0^{1/\eps}
h^2(\tilde{Z}^{\eps,x}_s)ds
+
\frac{Cq}{\eps}\int_0^{1/\eps}|\tilde{Z}^{\eps,x}_s|ds
\right]\right)^{\frac{1}{q(1+\delta)}}\\
&=
\left(\E\exp\left[
-\frac{p'q}{2\eps}\int_0^{1/\eps}
h^2(\tilde{Z}^{\eps,x}_s)ds
+
\frac{Cq}{\eps}\int_0^{1/\eps}|\tilde{Z}^{\eps,x}_s|ds
\right]\right)^{\frac{1}{q(1+\delta)}}\,.
\end{split}
\end{equation*}
Since $h(z)^2\ge h_0^2z^2$ (recall that $h(0)=0$!),
there exist $C(\delta)>0$, $C_1(\delta)$
such that $p'q h(z)^2/2-Cq|z|\geq C(\delta)z^2-C_1(\delta)$, 
and hence,
with $C_2(\delta)=C+C_1(\delta){\delta/p(1+\delta)}$ (all constants
here being positive and depending on $|||\tilde m|||, X$ only!),
\begin{equation}
\label{eq:tauold}
\begin{split}
\rho_1^\eps(x)&\le\exp\left[ \frac{C_2(\delta)}{\eps^2}+
\frac{\tilde{m}_0h(x)
%-h(\tilde{m}_0)x}{\eps} \right]
}{\eps} \right]
\left(\E
\left[p_0(\tilde{Z}^{\eps,x}_{1/\eps})\right]^{\frac{1+\delta}{\delta}}
\right)^{\frac{\delta}{1+\delta}}\\
&\quad\times
\left(\E\exp\left[-
\frac{C(\delta)}{\eps}\int_0^{1/\eps}|\tilde{Z}^{\eps,x}_s|^2ds
\right]\right)^{\frac{\delta}{q(1+\delta)}}\\
&\le \exp\left[\frac{C_3(\delta)}{\eps^2}\right]
\left(\E\exp\left[-
\frac{C(\delta)}{\eps}\int_0^{1/\eps}|\tilde{Z}^{\eps,x}_s|^2ds
\right]\right)^{\frac{\delta}{q(1+\delta)}}.
\end{split}
\end{equation}
It thus remains to estimate the last factor in 
the above right--hand side. Define 
$\tau=\inf\{t>0: |\tilde{Z}^{\eps,x}_s|<x/2\}$, and fix
$\eta>0$. We claim that for some $\eta>0$ small enough, it holds that
for some $C_\eta>0$, $x_0$ and all $|x|\geq x_0$,
\begin{equation}
\label{eq:taunew}
\P(\tau<\eta)\leq 
\exp\left[-\frac{C_\eta x^2}{\eps}\right]
\end{equation}
Assume (\ref{eq:taunew}), which will be proved below, and
 note that on the event $\{\tau\geq \eta\}$
we have that $\inf_{s\in (0,\eta]}
|\tilde{Z}^{\eps,x}_s|>x/2$.  We deduce from 
(\ref{eq:tauold}) 
\begin{equation}
\label{eq:tauold1}
%\begin{split}
\rho_1^\eps(x)%&
%\le\exp\left[ \frac{C}{\eps^2}+
%\frac{\tilde{m}_0h(x)-h(\tilde{m}_0)x}{\eps} \right]
%\left(\E\left[p_0(\tilde{Z}^{\eps,x}_{1/\eps})
%\right]^{\frac{1+\delta}{\delta}}
%\right)^{\frac{\delta}{1+\delta}}\\
%&\quad\times
%\left(\E\exp\left[-
%\frac{C(\delta)}{\eps}\int_0^{\eta}|\tilde{Z}^{\eps,x}_s|^2ds
%\right]\right)^{\frac{1}{1+\delta}}\\
\le
\exp\left[ \frac{C_3(\delta)}
{\eps^2}\right] %+\frac{\tilde{m}_0h(x)-h(\tilde{m}_0)x}{\eps} \right]
%\left(\E\left[p_0(\tilde{Z}^{\eps,x}_{1/\eps})\right]^{\frac{1+\delta}{\delta}}
%\right)^{\frac{\delta}{1+\delta}}\\
%&\quad
\times
\left(
\exp\left[-\frac{C_\eta x^2}{\eps}\right]+
\exp\left[-\frac{C(\delta)x^2\eta\delta}{4q(1+\delta)\eps}\right]
\right)\,,
%^{\frac{\delta}{q(1+\delta)}},
%\end{split}
\end{equation}
from which one concludes easily the bound
\begin{equation}
\label{eq:tauub}
\rho_1^\eps(x)\le\exp\left[ \frac{C_4(\delta)}{\eps^2}-
\frac{C x^2}{\eps}
\right]\,,
\end{equation}
for some constants $C_4(\delta)$ 
and $C$ depending on $\delta,|||\tilde m|||,X$ only.

On the other hand,
define the event
$$A_C=\{\sup_{t\in(0,1/\eps)} \sqrt{\eps} |\tilde{W}_t| \leq
C \}\,.$$
Then
there exists a constant $C_3>0$ depending on $C$ such that
$\P(A_C)\geq  C_3$. Note that on the
event $A_C$, because $h'(\cdot)>0$ and $h,b$ are Lipschitz,
Gronwall's inequality implies that
$\sup_{|x|\leq 1, s\leq 1/\eps}
|\tilde Z^{\eps,x}_s|\leq C'$ for some constant 
$C'$ depending on $C,\tilde m, X$ only.
Thus, on the event $A_C$,
%\be
%\label{eq-picard2}
%\rho^\eps_1(x):=
%e^{-F(x,\tilde m_0)/\eps}
$$|I_\eps(\tilde Z^{\eps,x}_{1/\eps}, 
0)+\int_0^{1/\eps} g_1(\tilde
Z_s^{\eps,x},\tilde m_s) ds+
\frac1{\eps} \int_0^{1/\eps} 
g_2(\tilde Z_s^{\eps,x},\tilde m_s) ds|
\leq \frac{C_4}{\eps^2}\,,
$$
where $C_4$ depends only on $\tilde m,X$ and 
the constants in Assumptions (A--1)-(A--4).
Hence, c.f. (\ref{eq-picard2}),
there exists a constant $C_2$ (again, depending on
the same quantities only)  such that uniformly in
$|x|<1$,
\begin{equation}
\label{eq:taulb}
\rho_1^\eps(x)\ge\exp\left[- \frac{C_2}{\eps^2}
\right]\,.
\end{equation}
(\ref{eq:taulb}) and (\ref{eq:tauub}) complete the proof of the 
lemma, once we prove (\ref{eq:taunew}).

Toward this end, assume without loss of generality that
$x>0$, 
and set $\hat h=2\sup_{y>0} h'(y)$.
Using the It\^o formula,
one has
\begin{equation}
\label{laila}
\tilde{Z}_t^{\eps,x} e^{\hat h t}=
x+\int_0^t
\left(\hat h \tilde{Z}_s^{\eps,x}-{h}(\tilde{Z}_s^{\eps,x})
+\tilde m_s
h'(\tilde{Z}_s^{\eps,x})-\eps b(\tilde{Z}_s^{\eps,x})\right)
e^{\hat h s}ds+\sqrt{\eps}\int_0^t e^{\hat h s} d\tilde{W}_s\,.
\end{equation}
Hence, denoting 
  $C_3=|||\tilde m||| \sup_xh'(x)
$, it follows that the event $\{\tau<\eta\}$ is
contained in the event
$$\{\sup_{t\in (0,\eta)}|\sqrt{\eps}\int_0^t
e^{\hat h s}
 d\tilde{W}_s|\geq x-C_3\frac{e^{\hat h \eta}-1}{\hat h}
-\frac{xe^{\hat h \eta}}{2}\}
\subset
\{\sup_{t\in (0,\eta)}|\sqrt{\eps}\int_0^t
e^{\hat h s}
 d\tilde{W}_s|
\geq \frac{x}{4}\}=:B\,,$$
if one choses $\eta$ small enough and $x$ large enough.
%Using
%that the process $Z_t:=\int_0^t \exp(\bar h(s-t)) d\tilde{W}_s$
%has the same law as an Orenstein-Uhlenbeck process,
We have that
$$\P(B)\leq
%\P\left(\sup_{t\in (0,\eta)}|Z_t|
%\geq \frac{e^{\bar h \eta}x}{4\sqrt{\eps}}\right)
%\leq
4 \exp\left(-\frac{Cx^2}{\eps}\right)\,,$$
for some constant $C$, which completes the proof of
(\ref{eq:taunew}).
$\Box$

\noindent
{\it Proof of Lemma \ref{lem-coarse2}:}
We only prove (\ref{eq-trajcut}), the proof of (\ref{eq-kijung}) being similar.
All we need to show is that for all $\eps\le\eps_0$, $|x|\le M_1/\sqrt{\eps}$,
and some $M_2$,
\begin{equation}\label{coarse2}
\E\left[L_\eps(x,1/\eps){\bf1}_{\{\|\tilde{Z}^{\eps,x}\|> M_2/\eps\}}\right]
\le
 \E\left[L_\eps(x,1/\eps){\bf1}_{\{\|\tilde{Z}^{\eps,x}\|\le M_2/\eps\}}\right].
\end{equation}
We first bound the left hand side of \eqref{coarse2} for $\eps\le1$. 
Recall the function $H$ introduced in (\ref{minnts}),
and apply
It\^o's formula to develop ${H}(\tilde{Z}^{\eps,x}_t)$
between $t=0$ and $t=1/\eps$, obtaining
\begin{equation*}
\begin{split}
\log L_\eps(x,1/\eps)-{H}(x)/\eps&=
-\frac{1}{2\eps}\int_0^{1/\eps}|h(\tilde{Z}^{\eps,x}_t)-
h(\tilde{m}_t)|^2dt
-\frac{1}{2\eps}\int_0^{1/\eps}
|{h}(\tilde{Z}^{\eps,x}_t)|^2dt\\
&+\frac{1}{\sqrt{\eps}}\int_0^{1/\eps}
{h}(\tilde{Z}^{\eps,x}_t)d\tilde{W}_t
+\int_0^{1/\eps}g_{3,\eps}(\tilde{Z}^{\eps,x}_t,\tilde{m}_t)dt
+\log p_0(\tilde{Z}^{\eps,x}_{1/\eps}),
\end{split}
\end{equation*}
where
\[
g_{3,\eps}(z,m)=g_1(z,m)-b(z){h}(z)+\frac{1}{2}h'(z)
+\frac{1}{2\eps}m^2(h'(z))^2.\]
Note that $\log p_0(\cdot)$ is bounded above, and 
\[
|g_{3,\eps}(z,\tilde{m}_t)|\le C\left(\frac{1}{\eps}+|z|\right).
\]
Now since for any $p>1$,
\[
\E\left[\exp\left(-\frac{p^2}{2\eps} \int_0^{1/\eps}|
{h}(\tilde{Z}^{\eps,x}_t)|^2dt
+\frac{p}{\sqrt{\eps}}\int_0^{1/\eps}{h}(\tilde{Z}^{\eps,x}_t)
d\tilde{W}_t
\right)\right]=1,\]
it follows from 
H\"older's inequality, that for any $q>p>1$ satisfying
$1/p+1/q=1$,
\begin{equation}\label{holdold}
\begin{split}
&e^{-{H}(x)/\eps}
\E\left[L_\eps(x,1/\eps){\bf1}_{\{\|\tilde{Z}^{\eps,x}\|> M_2/\eps\}}\right]
\le \\
&\left(\E\left[{\bf1}_{\{\|\tilde{Z}^{\eps,x}\|> M_2/\eps\}}
\exp\left( C\int_0^{1/\eps}(\frac{1}{\eps}+|\tilde{Z}^{\eps,x}_t|^2)dt \right)
\right.\right.
\\
&\left.\left.\times
\exp\left(-\frac{q}{2\eps}\int_0^{1/\eps}
|h(\tilde{Z}^{\eps,x}_t)-h(\tilde{m}_t)|^2dt
+\frac{p}{2\eps}\int_0^{1/\eps}|h(\tilde{Z}^{\eps,x}_t)|^2dt
\right)\right]\right)^{1/q}
%\\
%&\le
%\left(\E\left[{\bf1}_{\{\|\tilde{Z}^{\eps,x}\|> M_2/\eps\}}
%\exp\left(\frac{C}{\eps^2}
%\right)
%\right]
%\right)^{1/q},
\end{split}
\end{equation}
where $C>0$. But, note that due to $h'\geq h_0$, there
exists a constant $C$ depending on $|||\tilde m|||$
such that
$$ \sup_{z\in\reels, |m|\leq
|||\tilde m|||}
|z|^2-\frac{q}{2}|h(z)-h(m)|^2+
\frac{p}{2}|h(z)|^2\leq C\,.$$
Substituting in (\ref{holdold}),
one deduces that
\begin{equation}\label{hold}
e^{-{H}(x)/\eps}
\E\left[L_\eps(x,1/\eps){\bf1}_{\{\|\tilde{Z}^{\eps,x}\|> M_2/\eps\}}\right]
\le 
\left(\E\left[{\bf1}_{\{\|\tilde{Z}^{\eps,x}\|> M_2/\eps\}}
\exp\left(\frac{C}{\eps^2}
\right)
\right]
\right)^{1/q} \,.
\end{equation}
(Recall that the value of $C$ may change from line to line!).

We prove below that
%now deduce from the last inequality that, 
provided 
$M_2$ is large enough, there exists a $c>0$ such that
\begin{equation}\label{upnew}
%\E\left[L_\eps(x,1/\eps){\bf1}_{\{\|\tilde{Z}^{\eps,x}\|> 
%M_2/\eps\}}\right]
\E\left[{\bf1}_{\{\|\tilde{Z}^{\eps,x}\|> M_2/\eps\}}\right]
\le \exp\left(-\frac{c}{\eps^3}\right).
\end{equation}
Combined with (\ref{hold}), this implies that
uniformly in $|x|\leq M_1/\sqrt{\eps}$, 
\begin{equation}\label{up}
\E\left[L_\eps(x,1/\eps){\bf1}_{\{\|\tilde{Z}^{\eps,x}\|> 
M_2/\eps\}}\right]
\le \exp\left(-\frac{c}{\eps^3}\right).
\end{equation}
To see (\ref{upnew}), let $H=\sup|h'|$, define $\theta_0=0$
and 
$$\tau_i=\inf\{t>\theta_{i-1}: |\tilde{Z}_t^{\eps,x}|>\frac{M_2}{2\eps}\},
\quad
\theta_i=\inf\{t>\tau_{i}:
 |\tilde{Z}_t^{\eps,x}|<\frac{M_2}{4\eps}\}\,.
$$
Setting $f(z,m)=-h(z)+mh'(z)-\eps b(z)$, we have that
for $|z|\in [M_2/4\eps,M_2/\eps]$, $t\leq 1/\eps$
and $\eps$ small enough, it holds that
$ h_0M_2/8\eps\leq |f(z,\tilde m_t)|\le 2HM_2/\eps$ and
$\mbox{\rm sign}f(z,\tilde m_t)=-\mbox{\rm sign}(z)$.
Then, choosing $\eta=(16H)^{-1}$, for each $i$ it holds that
\begin{equation}
\label{iraq1}
\begin{split}
\P\left(\theta_i-\tau_i<\eta,\, 
\sup_{t\in[\tau_i,\theta_i]}|\tilde{Z}^{\eps,x}_t|<M_2/\eps
\right)
&\leq 
\P\left(\sqrt{\eps}\sup_{0\leq t\leq \eta} |W_t|
\geq \frac{M_2}{4\eps}-\frac{2H\eta M_2}{\eps}\right)
\\
&\leq
\P\left(\sqrt{\eps}\sup_{0\leq t\leq \eta} |W_t|
\geq \frac{M_2}{8\eps}\right)\\
&\leq \exp\left(-\frac{cM_2^2}{\eps^3 \eta }\right)\,.
\end{split}
\end{equation}
Similarly
\begin{equation}
\label{iraq11}
\begin{split}
\P\left(\theta_i-\tau_i\geq \eta,\, 
|\tilde{Z}^{\eps,x}_{\tau_i+\eta}|\geq M_2/2\eps
\right)
&\leq 
\P\left(\sqrt{\eps} W_\eta
\geq \frac{h_0 M_2\eta}{8\eps}\right)\\
&\leq
 \exp\left(-\frac{cM_2^2\eta}{\eps^3}\right)\,,
\end{split}
\end{equation}
and
\begin{equation}
\label{iraq2}
\begin{split}
\P\left(\sup_{t\in [\tau_i,(\tau_i+\eta)\wedge \theta_i]}
|\tilde{Z}_t^{\eps,x}|>M_2/\eps\right)
&\leq 
\P\left(\sqrt{\eps}\sup_{0\leq t\leq \eta} |W_t|
\geq \frac{M_2}{2\eps}\right)\\
&\leq \exp\left(-\frac{cM_2^2}{\eps^3\eta}\right)\,.
\end{split}
\end{equation}
Hence, using (\ref{iraq1}), (\ref{iraq11}) and (\ref{iraq2}),
$$\E\left[{\bf1}_{\{\|\tilde{Z}^{\eps,x}\|> M_2/\eps\}}\right]
\leq 
\frac{1}{\eps \eta}
\left(
\exp\left(-\frac{cM_2^2}{\eps^3\eta}\right)+
 \exp\left(-\frac{cM_2^2}{\eps^3}\right)
 +\exp\left(-\frac{cM_2^2\eta}{\eps^3 \eta }\right)
\right)\,,
$$
completing the proof of (\ref{upnew}).

%
%Let
%\[
%\tau_\eps:=\inf\{t\ge0;\, |\tilde{Z}_t^{\eps,x}|>\frac{M_2}{2\eps}\}.
%\]
%The next inclusion (here $\mu$ denotes the Lebesgue measure)
%follows by comparing the process $\tilde{Z}$ with Ornstein--Uhlenbeck~:
%\begin{equation*}
%\{\|\tilde{Z}^{\eps,x}\|> M_2/\eps\}\subset A_\eps\cup B_\eps\cup C_\eps,
%\end{equation*}
%where
%\[
%A_\eps=
%\{\mu\{0\le t\le\frac{1}{\eps};\, 
%|\tilde{Z}_t^{\eps,x}|>\frac{M_2}{4\eps}-M_3\}>1\},
%\]
%\[
%B_\eps=
%\{\sup_{t\le1/\eps-1}|\tilde{Z}_t^{\eps,x}|>\frac{M_2}{2\eps};
%\sup_{\tau_\eps\le t\le\tau_\eps+1}|\sqrt{\eps}(\tilde{W}_t
%-\tilde{W}_{\tau_\eps})|>\frac{M_2}{4\eps}-M_3\},
%\]
%and
%\[
%C_\eps=
%\{\sup_{1/\eps-1\le t\le1/\eps}|\sqrt{\eps}(\tilde{W}_t-\tilde{W}_{1/\eps-1})|
%>\frac{M_2}{2\eps}\}.
%\]
%On $A_\eps$, the right hand side of \eqref{hold} is bounded from above by
%$\exp(-c/\eps^3)$. 
%
%\[
%B_\eps\subset
%\{\sup_{\tau_\eps\le t\le\tau_\eps+1}|\sqrt{\eps}(\tilde{W}_t
%-\tilde{W}_{\tau_\eps})|>\frac{M_2}{4\eps}\},
%\]
%and the probability of this event is bounded from above by $\exp(-c/\eps^3)$,
%as well as the probability of the event $C_\eps$. The upper bound \eqref{up} 
%is established.

We now turn to the lower bound of the right hand side of \eqref{coarse2}.
Let, with $M'_1=M_1+1$,
\[
\eps_0=1\wedge \left(\frac{M_2}{M'_1}\right)^2.
\]
For $\eps\le\eps_0$,
\[
\{\|\tilde{Z}^{\eps,x}\|\le M'_1/\sqrt{\eps}\}
\subset
\{\|\tilde{Z}^{\eps,x}\|\le M_2/\eps\},
\]
so that for some $c'>0$
\begin{equation}\label{low}
\begin{split}
\E\left[L_\eps(x,1/\eps){\bf1}_{\{\|\tilde{Z}^{\eps,x}\|\le M_2/\eps\}}\right]
&\ge
\E\left[L_\eps(x,1/\eps){\bf1}_{\{\|\tilde{Z}^{\eps,x}\|\le M_1'/\sqrt{\eps}
\}}\right]\\
&\ge\exp\left(-\frac{c'}{\eps^{5/2}}\right)\P\left(\|\tilde{Z}^{\eps,x}\|\le M'_1/\sqrt{\eps}  \right)
\end{split}
\end{equation}
Finally \eqref{coarse2} follows from \eqref{up}, \eqref{low} and
the estimate
\[
\P\left(\|\tilde{Z}^{\eps,x}\|\le M'_1/\sqrt{\eps}  \right)\ge
\P\left(\sqrt{\eps}\|\tilde{W}\|\le C\right)\ge c''>0
.\]
\hfill $\Box$

\noindent
{\it Proof of Lemma \ref{lem-comp}:}
Note first that because of $(A-4)$, there exists a constant
$\kappa=\kappa(|||\tilde
m|||)$ such that for all $z\not\in [-\kappa,\kappa]$,
all $\eps<1/\kappa$, all $|m|\leq |||\tilde
m|||$,
and all $z'$,
$$\Delta(z,z',m)=-h(z)+h(z')+m[h'(z)-h'(z')]-\eps[
b(z)-b(z')]$$
satisfies $\sign(\Delta(z,z',m))=\sign(z'-z)$, while
$|\Delta(z,z',m)|\geq h_0 |z-z'|/2\,.$

Assume, w.l.o.g., that $x<y$.
Fix, for $\delta$ given,
a smooth, even, non-negative  function $c(z)$ such that
$c(|z|)$ is non-increasing,
$c(z)=\sqrt{\delta}$ for $|z|\leq \kappa$ and
$c(z)=0$ for $|z|>2\kappa$, with $||c'||\leq 10 \sqrt{\delta}$.
Define next the diffusions
\begin{eqnarray*}
d\xi^1_s&=&[-h(\xi_s^1)+\tilde m_s h'(\xi^1_s)
-\eps b(\xi_s^1)+ c(\xi_s^1){\bf 1}_{\{\tau>s
\}}]ds + \sqrt{\eps} dB_s\,
\quad \xi_0^1=x\,,\\
d\xi^2_s&=&[-h(\xi_s^2)+\tilde m_s h'(\xi^2_s)
-\eps b(\xi_s^2)]dt + \sqrt{\eps} dB_s\,
\quad \xi_0^2=y\,,
\end{eqnarray*}
where $B$ is a Brownian motion independent of the process
$\tilde m$, and $\tau=\min\{t: \xi_t^1=\xi_t^2\}\wedge 1/\eps$.
Note that $\xi^2$ coincides in distribution with 
$\tilde Z^{\eps,y}$, whereas the law of $\xi^1$
is %mutually 
absolutely continuous with respect to the law of $\tilde Z^{\eps,x}$
with Radon-Nykodim derivative given by
\begin{eqnarray}
\label{eq-lambda}
\Lambda&=&
\exp\left(\frac1{\eps}\int_0^\tau 
c(\xi_s^1)d\xi_s^1-\frac1{2\eps} \int_0^\tau c^2(\xi_s^1) ds
-\frac1{\eps}\int_0^\tau c(\xi_s^1) g(s,\xi_s^1)ds
\right)\\
&=&
\exp\left(\frac1{\eps}[\bar c(\xi_\tau^1)-\bar c(\xi_0^1)]
-\frac1{2\eps} \int_0^\tau c^2(\xi_s^1) ds
-\frac1{\eps}\int_0^\tau c(\xi_s^1) g(s,\xi_s^1)ds
-\frac1{2}\int_0^\tau c'(\xi_s^1) ds
\right)
\,,\nonumber \end{eqnarray}
where $g(s,z)=-h(z)+\tilde m_sh'(z)-\eps b(z)$ and
$\bar c(z)=\int_0^z c(y) dy$.

Next, note that with $\zeta_s=\xi_s^1-\xi_s^2$, and using that 
$x<y$, it holds that $\zeta_s\leq 0$ for all $s$, while by 
definition $|\zeta_0|\leq \delta$. 
Hence,
by the definition of $c(\cdot)$ and of $\kappa$, it holds that
for all $\delta<\delta_1(\kappa,|||m|||)$,
$$ d\zeta_s/ds\geq -\frac{h_0\zeta_s}{2}+\frac{c(\xi_s^1){\bf 1}_{s<\tau}}{2},$$
from which one concludes that $\zeta_s\geq -\delta e^{-hs/2}$.
In particular, this implies that
for all such $\delta$,
$$\int_0^\tau c(\xi_s^1){\bf 1}_{\{\tau>s\}} ds
=\int_0^\tau c(\xi_s^1) ds
\leq C \delta$$
for some constant $C=C(\kappa,|||\tilde m|||)$.
Since $c(z)=0$ for $|z|>2\kappa$, and since
$|g(s,z)|$ is bounded uniformly in $s\leq 1/\eps$ and
$|z|\leq 2\kappa$ (by a bound that depends only on 
$|||\tilde m|||$), the last inequality implies that
$$|\int_0^\tau c(\xi_s^1) g(s,\xi_s^1)ds|
\leq C\delta$$
again, for some constant $C$ depending on $\kappa,|||\tilde m|||$
only. 
Finally, note that
$$
\int_0^\tau c^2(\xi_s^1) ds\leq
\sqrt{\delta}\int_0^\tau c(\xi_s^1) ds\leq
C\delta^{3/2}\,,$$
and that $|\bar c(z)|\leq 2\kappa \sqrt{\delta}$.
Substituting back into 
(\ref{eq-lambda}), and recalling that
$\kappa=\kappa(|||\tilde m|||)$,
one concludes the existence of a constant
$C_2=C_2(|||\tilde m|||)$ such that for all $\delta<\delta_1$,
\begin{equation}
\label{eq-lambda1}
e^{-C_2\sqrt{\delta}/\eps}\leq \Lambda
\leq
e^{C_2\sqrt{\delta}/\eps}
\,.
\end{equation}
Therefore, with $\E_B$ denoting expectation with respect to 
$B_\cdot$, and using the bound on $\Lambda$ in the second
inequality,
and the Lipschitz property of $g_1,g_2$ together with
the exponential decay of $\zeta_s$ in the third, that
for all $t>1/2\eps$, and ommiting the dependence on
$\theta^{1/\eps-t}\tilde m$ everywhere,
\begin{eqnarray}
\E L_\eps(x,t)&\leq& 
2\E L_\eps(x,t)
{\bf 1}_{\{||\tilde Z^{\eps,x}||<M_2/\eps\}}
\nonumber \\
&=&2
\E_B\left(
{\bf 1}_{\{||\xi^1||<M_2/\eps\}}
\Lambda^{-1}  
\exp\left(
I_\eps(\xi_{t}^1,0)+
\int_0^{t} \left(g_1(\xi_s^1,\tilde m_s)+
\frac1{\eps}g_2(\xi_s^1,\tilde m_s)\right)ds\right)\right)\nonumber \\
&\leq &2
\E_B\left(\left.
{\bf 1}_{\{||\xi^2||<(M_2+1)/\eps\}}
\right.\right.\nonumber\\
&& \;\left(
\exp\left(\frac{C_2\sqrt{\delta}}{\eps}+
I_\eps(\xi_{t}^2+\zeta_{t},0)+
\int_0^{t} \left(g_1(\xi_s^2+\zeta_s,\tilde m_s)+
\frac1{\eps}g_2(\xi_s^2+\zeta_s,\tilde m_s)\right)ds\right)\right)
\nonumber \\
&\leq &2
\E_B\left(\exp\left(\frac{C_3\sqrt{\delta}}{\eps}+
I_\eps(\xi_{t}^2,0)+
\int_0^{t} \left(g_1(\xi_s^2,\tilde m_s)+
\frac1{\eps}g_2(\xi_s^2,\tilde m_s)\right)ds\right)\right)
\nonumber \\
&= &2
\E\left(\exp\left(\frac{C_3\sqrt{\delta}}{\eps}+
I_\eps(\tilde Z_{t}^{\eps,y},0)+
\int_0^{t} \left(g_1(\tilde Z_s^{\eps,y},\tilde m_s)+
\frac1{\eps}g_2(\tilde Z_s^{\eps,y},\tilde m_s)\right)ds
\right)\right)
\nonumber\\
&=&2
\exp\left(\frac{C_3\sqrt{\delta}}{\eps}\right)
\E L_\eps(y,t)
\nonumber\\
&\leq & 
4\exp\left(\frac{C_3\sqrt{\delta}}{\eps}\right)
\E \left(
{\bf 1}_{\{||\tilde Z^{\eps,y}||<M_2/\eps\}}
L_\eps(y,t)\right)
\,,
\end{eqnarray}
yielding (\ref{eq-comp1}) 
for $x<y$ and $\delta<\delta_1$, with $g(\delta)=
C_3\sqrt{\delta}$. 
Further, the same computation gives
$$
4 \E \left( L_\eps(x,t)
{\bf 1}_{\{||\tilde Z^{\eps,x}||<M_2/\eps\}}\right)
\geq
\exp\left(\frac{-C_3\sqrt{\delta}}{\eps}\right)
\E \left(L_\eps(y,t)
{\bf 1}_{\{||\tilde Z^{\eps,y}||<M_2/\eps\}}
\right)\,,
$$
yielding, by exchanging the roles of $x$ and $y$,
 (\ref{eq-comp1}) for $x>y$ and $\delta<\delta_1$ with the
same $g(\delta)$. Finally, for $\delta>\delta_1$, iterate this 
procedure to obtain (\ref{eq-comp1})
with $g(\delta)=C_3\sqrt{\delta\wedge \delta_1}\lceil \delta/\delta_1
\rceil$.
Substituting $y=X_1$ gives then (\ref{eq-comp2}).
$\Box$

\noindent
{\it Proof of Lemma \ref{lem-local}:}
%The $x$ dependence part of the claim follows the same lines as
%in the proof of Lemma \ref{lem-comp}.
%{\bf Add details}. We concentrate next on the proof of the 
%dependence in $z$, beginning with the upper bound. Toward
Throughout the proof, we fix once and for all the sequence 
$T_\eps$. All constants $C_i$ used in the
proof may depend on the choice of the
sequence  but not explicitely on $\eps$.

We begin with the proof of (\ref{eq-110}).
Using Girsanov's theorem one finds that
with $\bar Z_t^{\eps,x}=x+\sqrt{\eps}\tilde{W}_t$,
\begin{eqnarray}
\label{aa}
&&\E\left[\bar L_\eps(x,T)
{\bf 1}_{\{|\tilde Z_T^{\eps,x}-z|<\delta\}}
{\bf 1}_{\{||\tilde Z^{\eps,x}||_T\leq M_3/\eps\}}
\right]\\
&&
=
\E\left[
{\bf 1}_{\{|\bar Z_T^{\eps,x}-z|<\delta\}}
{\bf 1}_{\{||\bar Z^{\eps,x}||_T\leq M_3/\eps\}}
\exp\left(
\frac{1}{\eps}
\int_0^{T}\left[-h(\bar Z_s^{\eps,x})+
\tilde m_s h'(\bar Z_s^{\eps,x}) -\eps b(\bar Z_s^{\eps,x})\right]
d\bar Z_s^{\eps,x}\right.\right.\nonumber\\
&&\left.\left.\quad
 -\frac{1}{\eps}
\int_0^{T} \left(
\frac{[h(\bar Z_s^{\eps,x})-h(\tilde m_s)]^2}{2}+
\frac{b^2(\bar Z_s^{\eps,x}) \eps^2}{2}+
\eps b(\bar Z_s^{\eps,x})h(\bar Z_s^{\eps,x})
-\eps h'(\bar Z_s^{\eps,x})b(\bar Z_s^{\eps,x})\tilde m_s
\right.\right.\right.\nonumber\\
&&\left.\quad\quad
\quad\quad
-\eps g_1(\bar Z_s^{\eps,x},\tilde m_s)
\Bigl)ds\Bigl)\right]\nonumber
\end{eqnarray}
We consider separately the  different terms in (\ref{aa}).
Note first that one may, exactly as
in the course of the proof of
Lemma \ref{lem-comp},  move from starting point $x$
to starting point $X_1$ in
the right hand side of (\ref{aa}), with the effect of
picking up a term bounded by
$\exp (C|x|/\eps)$ and
widening the allowed region where $\bar Z_T^{\eps,x}$ need to be,
namely for
all $T_\eps\geq T>1$,
the right hand side of (\ref{aa}) is bounded by
\begin{eqnarray}
\label{aa1}
&&\exp\left(
\frac{C_1+C_2|x|}{\eps}\right)
\E\left[
{\bf 1}_{\{|\bar Z_T^{\eps,X_1}-z|<\delta+|x|+|X_1|\}}
{\bf 1}_{\{||\bar Z^{\eps,X_1}||_T\leq |x|+(M_3+1)/\eps\}}\right.\\
&&\quad \quad\quad
\left.\exp\left(
\frac{1}{\eps}
\int_0^{T}\left[-h(\bar Z_t^{\eps,X_1})+
\tilde m_s h'(\bar Z_s^{\eps,X_1}) -\eps b(\bar Z_s^{\eps,X_1})\right]
d\bar Z_s^{\eps,X_1}\right)\right]\,. \nonumber
\end{eqnarray}
An integration by parts gives that
$$-\int_0^{T}h(\bar Z_t^{\eps,X_1})
d\bar Z_t^{\eps,X_1}=
-\bar{\cal{J}}(\bar Z_T^{\eps,X_1},X_1)-h(X_1)(\bar Z_T^{\eps,X_1}-X_1)+
\frac{\eps}{2} \int_0^T h'(\bar Z_t^{\eps,X_1})dt\,,
$$
and hence, on the event $
\{|\bar Z_T^{\eps,X_1}-z|<\delta+|x|+|X_1|\}$,
it holds that
\begin{equation}
\label{aa2}
-\int_0^{T}h(\bar Z_t^{\eps,X_1})d\bar Z_t^{\eps,X_1}\leq
-C(|z|-|x|-|X_1|-\delta)_+^2+C.
\end{equation}
Similarly, with $B(z)=\int_{X_1}^z b(x) dx$,
\begin{equation}
\label{aa3}
\int_0^T 
 b(\bar Z_s^{\eps,X_1})
d\bar Z_s^{\eps,X_1}=
B(\bar Z_T^{\eps,X_1})-\frac{\eps}2 
\int_0^T b'(\bar Z_s^{\eps,X_1})ds\leq
C(|z|^2+|x|^2+1)\,.
\end{equation}
Finally, rewrite
$$\int_0^{T}
\tilde m_s h'(\bar Z_s^{\eps,X_1})
d\bar Z_s^{\eps,X_1}=
X_1 \int_0^{T}
h'(\bar Z_s^{\eps,X_1})
d\bar Z_s^{\eps,X_1}+
\int_0^{T}
(\tilde m_s-X_1) h'(\bar Z_s^{\eps,X_1})
d\bar Z_s^{\eps,X_1}\,.$$
The first stochastic integral in the above expression is
handled exactly as in (\ref{aa3}), and substituting
in (\ref{aa1})
one concludes that the right hand side of
(\ref{aa}) is bounded by
\begin{eqnarray*}
&&
\exp\left(
\frac{C+C(|x|+|z|)-C(|z|-|x|)_+^2}{\eps}\right)
\E\left[
\exp\left(
\frac{1}{\eps}\int_0^{T}
(\tilde m_s-X_1) h'(\bar Z_s^{\eps,X_1})
d\bar Z_s^{\eps,X_1}\right)\right]
\\
&&\leq
\exp\left(
\frac{C+C(|x|+|z|)-C(|z|-|x|)_+^2}{\eps}
+\frac{1}{2\eps}\int_0^{T_\eps} C|\tilde m_s-X_1|^2ds\right)\\
&&\leq 
\exp\left(
\frac{C+C(|x|+|z|)-C(|z|-|x|)_+^2}{\eps}\right)\,,
\end{eqnarray*}
where in the last inequality 
we used the last part of Lemma \ref{lemma-m}.
This completes the proof of (\ref{eq-110}).

The proof of (\ref{eq-111}) proceeds along similar lines.
The starting point is the change of measure leading to 
(\ref{aa}). Define the function
$$\Psi_t=\left\{\begin{array}{ll}
x+2(X_1-x)t, & t\leq 1/2\\
X_1, & T-1/2>t\geq 1/2\,,\\
z+2(z-X_1)(t-T), & T\geq t\geq T-1/2\,.
\end{array}
\right.$$
Let $D$ denote the
event
$$D:=\{\sup_{t\leq T} |\bar Z_t^{\eps,x} -\Psi_t|<\sqrt{\eps}\}\,.$$
We will prove below that for $|x-X_1|\leq 1$, and $T<T_\eps$,
there exists
a constant $C$ independent of $T$ and $\eps$ such that
\begin{equation}
\label{boarding}
\P(D)
\geq e^{- \frac{C}{\eps}}\,.
\end{equation}
%Assuming that ,
We can clearly bound from below the right hand side 
of 
(\ref{aa}) by
\begin{eqnarray*}
&&\E\left[
{\bf 1}_{\{|\bar Z_T^{\eps,x}-z|<\delta\}}
{\bf 1}_{\{||\bar Z^{\eps,x}||_T\leq M_3/\eps\}}
{\bf 1}_{D}
\exp\left(
\frac{1}{\eps}
\int_0^{T}\left[-h(\bar Z_t^{\eps,x})+
\tilde m_s h'(\bar Z_s^{\eps,x}) -\eps b(\bar Z_s^{\eps,x})\right]
d\bar Z_s^{\eps,x}\right.\right.\\
&&\left.\left.\quad
 -\frac{1}{\eps}
\int_0^{T} \left(
\frac{[h(\bar Z_t^{\eps,x})-h(\tilde m_t)]^2}{2}+
\frac{b^2(\bar Z_t^{\eps,x}) \eps^2}{2}+
\eps b(\bar Z_s^{\eps,x})h(\bar Z_s^{\eps,x})
-\eps h'(\bar Z_s^{\eps,x})b(\bar Z_s^{\eps,x})\tilde m_s
\right.\right.\right.\nonumber\\
&&\left.\quad\quad
\quad\quad
-\eps g_1(\bar Z_s^{\eps,x},\tilde m_s)
\Bigl)ds\Bigl)\right]\,.
\end{eqnarray*}
We now assume that (\ref{boarding}) and $|z-X_1|\le1$ hold. Then
using the same integration by parts as in the proof of the 
upper bound, one concludes that
the right hand side 
of 
(\ref{aa}) is bounded from below by
\begin{equation}
\label{flight}
\E\left[
{\bf 1}_{D}
\exp\left(
\frac{-C}{\eps}+ 
\frac{1}{\eps}\int_0^T(\tilde m_s-X_1)
h'(\bar Z_s^{\eps, x})d\bar Z_s^{\eps,x}\right)\right] \,.
\end{equation}
But, since
$$\mbox{\rm Var}
\left(\int_0^T(\tilde m_s-X_1)
h'(\bar Z_s^{\eps, x})d\bar Z_s^{\eps,x}\right)\leq 
C\eps \,,$$
%
%{\bf I belive you were wrong with the $\eps^2$. Now I am not so sure how to finish.
% If $X$ denotes the above stochastic integral, we have
% $$\E\left[\exp(X/\eps)|D\right]\ge \exp(-c/\eps)\P[X>-c|D],$$
% and we need to estimate $\P(\{X>-c\}\cap D)$. Easy ? Hard ? }
%
one gets,
using Chebycheff's inequality, that
$$\P
\left[
\int_0^T(\tilde m_s-X_1)
h'(\bar Z_s^{\eps, x})d\bar Z_s^{\eps,x}<-c
\right]
\leq
\exp \left(-\frac{C_2c^2}{\eps}\right)\,.$$
Hence,
$$\P
\left[
\int_0^T(\tilde m_s-X_1)
h'(\bar Z_s^{\eps, x})d\bar Z_s^{\eps,x}<-c|D
\right]
\leq \frac{
\exp \left(-\frac{C_2c^2}{\eps}\right)}{\P(D)}
\leq \frac{1}{2}\,,$$
if $c$ is chosen large, where in the last inequality we used
(\ref{boarding}). In particular, it follows that
$$\E\left[
\exp\left(
\frac{1}{\eps}\int_0^T(\tilde m_s-X_1)
h'(\bar Z_s^{\eps, x})d\bar Z_s^{\eps,x}\right)\,|\,D\right]
\geq \exp\left(-\frac{C}{\eps}\right)
\,,
$$
for some $C>0$. Substituting back in (\ref{flight})
the required lower bound follows.

It thus only remains to prove (\ref{boarding}).
This however is immediate from a martingale argument: first, perform
the change of measure making $
S_t:=\bar Z_t^{\eps,x} -\Psi_t$ into a Brownian motion of variance
$\eps$. Then, for $1\leq T\leq T_\eps$,
$$\P(D)=
\E\left({\bf 1}_{\{\sup_{t\leq T} |S_t|\leq  \sqrt{\eps}\}}
\exp\left(-\frac{1}{\eps}
\int_0^{T} \dot{\Psi_t}dS_t-\frac1{2\eps}
\int_0^{T} \dot{\Psi_t}^2 dt\right)\right)\,.$$
Integrating by parts the stochastic integral, and using
that $\dot{\Psi}(t)=0$ for $t\in (1/2,T-1/2)$,
(\ref{boarding}) follows, which completes the proof of the lemma.
$\Box$

\noindent
{\it Proof of (\ref{exp-tight})}
We let $\eta>0$ as before.
Note first that by (\ref{eq:tauub}) and (\ref{eq:taulb}),
there is a constant $M$ depending on $|||\tilde m|||$
only such that
\begin{equation}
\label{outofrange}
\limsup_{\eps\to0}
\eps\log \int_{[-M/\sqrt{\eps},M/\sqrt{\eps}]^c}
q_1^\eps(x)dx =-\infty\,.
\end{equation}
We may and will in the sequel assume that $M=M_1$ where
$M_1$ is defined in Lemma \ref{lem-coarse1}, and
we use  $M_3$ and
$M_2$ as in Lemma \ref{lem-coarse2}. 

Next, 
set $\eps_4$ such that $\eps_4\log 2<\eta/8$ and
$\eps\log(2M_3/\eps\delta)\leq \eta/8$
for $\eps<\eps_4$.
Repeating the arguments in (\ref{eq-fin1}),
without using the compact set ${\cal K}_1$, one has
%by (\ref{eq-fin1night}),
for $\eps<\eps_4$ and $|x|\leq M_1/\sqrt{\eps}$,
\begin{eqnarray}
\label{eq-tightnight}
\eps\log \rho_1^\eps(x)&\leq&
-F(x,\tilde m_0)+
\eps\log \tilde J_\eps(x)
+\frac{\eta}{4}\quad\mbox{\tt 
as in (\ref{eq-fin1})}\nonumber\\
&\leq &
-F(x,\tilde m_0)+
\eps\log \sup_{|z|\leq M_3/\eps}
\hat J_{\eps,T}(x,z)
+\frac{\eta}{2}\quad\mbox{\tt 
by (\ref{eq-revmor2}) and (\ref{turkey1})}\nonumber \\
&\leq &
-F(x,\tilde m_0)+\frac{\eta}{2}+C_2-C_3(|z|-|x|)_+^2
+C_5(|x|+|z|)\nonumber \\
&&\quad 
+\eps\log
%\mbox{\tt by \ref{eq-rev10}}\nonumber\\
%&&\quad \cdot
\E \left[L_\eps(X_1,1/\eps-T,\theta^T \tilde m)
{\bf 1}_{\{||\tilde Z^{\eps,X_1}||_{1/\eps-T}\leq M_2/\eps\}}\right]\,.
\end{eqnarray}
A similar argument shows that for $|x-X_1|<1$,
and some constant $C_6$ depending only on $X,|||\tilde m|||$,
\begin{equation}
\label{eq-tightnight1}
\eps\log \rho_1^\eps(x)\geq
-F(X_1,X_1)-C_6
+\eps\log
\E \left[L_\eps(X_1,1/\eps-T,\theta^T \tilde m)
{\bf 1}_{\{||\tilde Z^{\eps,X_1}||_{1/\eps-T}\leq M_2/\eps\}}\right]\,.
\end{equation}
Fixing now an $L$, and using as in (\ref{eq-rev12}) the uniform
quadratic 
growth  of $F(x,m)$ as $|x|\to\infty$ and
$|m|<|||\tilde m|||$, one finds a compact set ${\cal K}^L$
such that
\begin{equation}
\label{eq-rev12night}
\sup_{|m|<|||\tilde m|||}
\sup_{x\in ({\cal K}^L)^c, z\in\reels}
\frac{C_2}{\eps}-\frac{C_3(|z|-|x|)_+^2}{\eps}
+\frac{C_5(|x|+|z|)}{\eps} -F(x,m)
\leq -F(X_1,X_1) -\frac{C_6+L}{\eps}\,,
\end{equation}
and hence, from (\ref{eq-tightnight}) and  (\ref{eq-tightnight1}),
for $x\in ({\cal K}^L)^c\cap [-M_1/\sqrt{\eps},M_1/\sqrt{\eps}]$,
\begin{equation}
\label{almost}
\eps\log \rho_1^\eps(x)\leq \inf_{|y-X_1|\leq1}
\eps\log \rho_1^\eps(y)-L\,.
\end{equation}
Hence,
\begin{eqnarray}
\limsup_{\eps\to 0}
\eps\log\int_{({\cal K}^L)^c} q_1^\eps(x)dx&
= &
\limsup_{\eps\to 0}
\eps\log\int_{({\cal K}^L)^c\cap[-M_1/\sqrt{\eps},M_1/\sqrt{\eps}]} 
q_1^\eps(x)dx
\quad \mbox{\tt by (\ref{outofrange})}
\nonumber
\\
&\leq&
\limsup_{\eps\to 0}\left[
\eps\log\int_{({\cal K}^L)^c\cap[-M_1/\sqrt{\eps},M_1/\sqrt{\eps}]} 
\rho_1^\eps(x)dx-
\eps\log\int_{[X_1-1,X_1+1]} 
\rho_1^\eps(x)dx\right]\nonumber\\
\mbox{\tt by (\ref{almost})}
&\leq &
\limsup_{\eps\to 0}\left[
\eps\log\left(\frac{2M_1}{\sqrt{\eps}}\right)+
\inf_{|y-X_1|\leq 1}\eps\log \rho_1^\eps(y)-L
-\inf_{|y-X_1|\leq 1}\eps\log \rho_1^\eps(y)
\eps\log 2
\right]\nonumber \\
&\leq & -L \,.
\end{eqnarray}
This completes the proof.
$\Box$.

\noindent
{\bf Acknowledgement} We thank Ki-Jung Lee 
for a careful 
reading of a preliminary version
of this paper. We also thank an anonymous referee for a detailed
reading of the paper, and many useful and important comments.

\section*{Appendix: derivation of
(\ref{eq-picard1})}
\setcounter{equation}{0}
We first recall Picard's theorem, \cite[Proposition 4.2]{picard}:
under the assumptions of the current paper and with the same notations,
a version of the conditional unnormalized density is given by
\begin{equation}
\label{eq-pic1}
%\begin{split}
\tilde q(1,x)=\exp\left\{\frac{1}{2\eps^2}\int_0^1 h^2(\bar m_s))ds -\frac{1}{\eps} 
 F(x,\tilde m_0)\right\}
 \tilde \E'\left[\exp \rho_1^{y,x}\right]
\end{equation}
where
\begin{eqnarray*}
%\bar F(s,x)&=&F(x,\bar m_s)\,,\\
\rho_1^{x,y}&=&
\log p_0(\bar X_{1}^x)+\frac{1}{\eps}F(\bar X_1^x,0)
-\frac{1}{\eps}\int_0^1 h(\bar m_s) d\bar X_s^x
-\frac{1}{\eps}\int_0^1 h(\bar X_s^x) b(\bar m_s) ds
+\frac{1}{\eps} \int_0^1 \bar m_s  h'(\bar X_s^x) d\bar X_s^x
\\
&& 
+\frac{1}{2\eps} \int_0^1 \bar m_s  h''(\bar X_s^x) ds
+\frac{1}{\eps}\int_0^1\left[b(\bar X_s^x)(h(\bar X_s^x)-h(\bar m_s))
-\frac12 h'(\bar X_s^x)-\eps b'(\bar X_s^x)\right]ds\\
d\bar X_s^x&=&
- \frac{1}{\eps} (h(\bar X_s^x)-h(\bar m_s)) ds
-b(\bar X_s^x) ds + dW_s\,, \quad \bar X_0^x=x\,,
\end{eqnarray*}
$W_\cdot$ is a Brownian motion, and $\tilde \E'$ denotes expectation
with respect to this Brownian motion. Performing a time change $t\mapsto 
\eps t$ and setting $\tilde W_t= \frac{1}{\sqrt{\eps}} W_{\eps t}$,
we have that $\tilde W_t$ is again a standard Brownian motion and
with $\bar X_t^{\eps,x}=\bar X_{\eps t}^x$,
\begin{eqnarray*}
\rho_1^{x,y}&=&
\log p_0(\bar X_{1/\eps}^{\eps, x})+\frac{1}{\eps} F(\bar X_{1/\eps}^{\eps,x},0)
-\frac{1}{\eps}\int_0^{1/\eps} h(\tilde m_s) d\bar X_s^{\eps,x}
-\int_0^{1/\eps} h(\bar X_s^{\eps,x}) b(\tilde m_s) ds
\\
&&+ \frac{1}{\eps} \int_0^{1/\eps} \tilde m_s  h'(\bar X_s^{\eps,x}) d\bar  X_s^{\eps,x}
+ \frac{1}{2} \int_0^{1/\eps} \tilde m_s  h''(\bar X_s^{\eps,x}) ds
\\
&&+\int_0^{1/\eps}\left[b(\bar X_s^{\eps,x})(h(\bar X_s^{\eps,x})-h(\tilde m_s))
-\frac12 h'(\bar X_s^{\eps,x})-\eps b'(\bar X_s^{\eps,x})\right]ds\\
d\bar X_s^{\eps,x}&=&
-  (h(\bar X_s^{\eps,x})-h(\tilde m_s)) ds
-\eps b(\bar X_s^{\eps,x}) ds +\sqrt{\eps} d\tilde W_s\,, \quad \bar X_0^{\eps,x}=x\,,
\end{eqnarray*}
and 
\begin{equation}
\label{eq-pic2}
%\begin{split}
\tilde q(1,x)=\exp\left\{\frac{1}{2\eps}\int_0^{1/\eps} h^2(\tilde m_s))ds -\frac{1}{\eps} 
F(x,\tilde m_0)\right\}
 \tilde \E\left[\exp \rho_1^{y,x}\right]\,,
\end{equation}
where the expectation now is with respect to the Brownian motion $\tilde W_t$.

Observe next that, by Girsanov's theorem, the law of  the process $\bar X_t^{\eps,x}$
is absolutely continuous with respect to that of the process $\tilde Z_t^{\eps,x}$, with 
Radon-Nykodym derivative given by
\begin{equation}
\begin{split}
&e^{\Lambda}=
\exp\left[\frac{1}{{\eps}}
\int_0^{1/\eps}
[h(\tilde m_s)-\tilde m_s h'(\tilde Z_s^{\eps,x})]
d\tilde Z_s^{\eps,x}-
\frac{1}{2\eps} \int_0^{1/\eps}[h(\tilde Z_s^{\eps,x})-h(\tilde m_s)+\eps b(\tilde Z_s^{\eps,x})
]^2 ds\right.\\
&\quad \quad \quad \left. 
+\frac{1}{2\eps} \int_0^{1/\eps}[h(\tilde Z_s^{\eps,x})-\tilde m_s h'(\tilde Z_s^{\eps,x})
+\eps b(\tilde Z_s^{\eps,x})
]^2 ds\right]\,.
\end{split}
\end{equation}
Hence, with $\E$ denoting expectations with respect to the Brownian motion
$\tilde W_t$ appearing in the definition of
$\tilde Z_t^{\eps,x}$, (\ref{eq-pic2}) transforms to
$$\tilde q(1,x)=\exp\left\{\frac{1}{2\eps}\int_0^{1/\eps} h^2(\tilde m_s))ds -\frac{1}{\eps} 
F(x,\tilde m_0)\right\}\E \exp[\Lambda_1(x)]\,,$$
where 
\begin{eqnarray*}
\Lambda_1(x)&=&
\Lambda +
\log p_0(\tilde Z_{1/\eps}^{\eps, x})+\frac{1}{\eps} F(\tilde Z_{1/\eps}^{\eps,x},0)
-\frac{1}{\eps}\int_0^{1/\eps} h(\tilde m_s) d\tilde Z_s^{\eps,x}
-\int_0^{1/\eps} h(\tilde Z^{\eps,x}) b(\tilde m_s) ds
\\
&&
+\frac{1}{\eps}
 \int_0^{1/\eps} \tilde m_s  h'(\tilde Z_s^{\eps,x}) d\tilde Z_s^{\eps,x}
+\frac12 \int_0^{1/\eps} \tilde m_s h''(\tilde Z_s^{\eps,x}) ds\\
&&
+\int_0^{1/\eps}\left[b(\tilde Z_s^{\eps,x})(h(\tilde Z_s^{\eps,x})-h(\tilde m_s))
-\frac12 h'(\tilde Z_s^{\eps,x})-\eps b'(\tilde Z_s^{\eps,x})\right]ds\\
&=&
\log p_0(\tilde Z_{1/\eps}^{\eps, x})+\frac{1}{\eps} F(\tilde Z_{1/\eps}^{\eps,x},0)
+\int_0^{1/\eps} g_1(\tilde Z_s^{1/\eps},\tilde m_s)ds+
\frac{1}{\eps}\int_0^{1/\eps} g_2(\tilde Z_s^{1/\eps},\tilde m_s)ds\,.
\end{eqnarray*}
Since  $\int_0^{1/\eps} h^2(\tilde m_s) ds$ does not depend on $x$,
taking 
$$\rho_1^\eps(x)=
\tilde  q(1,x)\exp\left\{-\frac{1}{2\eps}\int_0^{1/\eps} h^2(\tilde m_s))ds\right\}\,,
$$
gives a version of the unnormalized conditional  density that coincides
with (\ref{eq-picard2}).
$\Box$.

\bibliographystyle{plain}

\end{document}